\documentclass[12pt,reqno]{amsart}

\usepackage{geometry}
\usepackage{verbatim}
\usepackage{curves}
\usepackage{longtable}
\usepackage{color}
\usepackage{epic}
\usepackage{amssymb}
\usepackage{amsthm}
\usepackage{amscd}
\usepackage{amsmath}
\usepackage{graphics}
\allowdisplaybreaks
\begin{document}

\newtheorem{main}{Theorem}

\newtheorem{Theorem}{Theorem}
\newtheorem{theorem}{Theorem}[section]
\newtheorem{corollary}[theorem]{Corollary}
\newtheorem{proposition}[theorem]{Proposition}
\newtheorem{lemma}[theorem]{Lemma}

\newtheorem{thmletter}{Theorem}
\def\thethmletter{\Alph{thmletter}}

\theoremstyle{definition}
\newtheorem{definition}[theorem]{Definition}
\newtheorem{example}[theorem]{Example}
\newtheorem{remark}[theorem]{Remark}
\newtheorem*{remark*}{Remark}
\newtheorem{notation}{Notation}[section]

\newcommand{\noqed}{\def\qedsymbol{}}
\def\mtline#1{\hbox to#1{\hrulefill}} 
\def\bbk{\bigbreak}
\let\bbk\relax
\def\Gam{\Gamma}
\def\Sig{\Sigma}
\def\lds{\dots}
\def\blt{\mbox{$\begin{picture}(6,6)\put(3,3){\circle*{6}}\end{picture}$}}
\def\noi{\noindent}
\let\noi\relax
\def\wtit{\widetilde}
\def\ub{\underline}
\def\Raw{\Rightarrow}
\def\lraw{\longrightarrow}
\def\ul{\underline}

\renewcommand{\qedsymbol}{Q.E.D.}

\newcommand{\zz}{\phantom{00}}
\newcommand{\TP}{\!{}+{}\!}
\newcommand{\TM}{\!{}-{}\!}
\newcommand{\TE}{\!{}={}\!}

\title[Normal Surface singularities of small degrees]{Normal Surface singularities of small degrees}

\author[S. Yau]{Stephen S.-T. Yau}
\address{ 
Beijing Institute of Mathematical Sciences and Applications (BIMSA),   Beijing 101400\\ P. R. of China;
Department of Mathematical Sciences,
	Tsinghua University,
	 Beijing 100084\\ P. R. of China} \curraddr{}
\email{yau@uic.edu}
%\thanks{Corresponding author: Stephen S.-T. Yau}
\author[H. Zuo]{Hao Zuo}
\address{Department of Mathematical Sciences\\ Tsinghua University\\
 Beijing 100084\\ P. R. of China} \curraddr{}
\email{zuoh22@mails.tsinghua.edu.cn}
%\thanks{Qiwei Zhu  is co-first author of this paper}
\author[H. Zuo]{Huaiqing Zuo}
\address{Department of Mathematical Sciences\\ Tsinghua University\\
 Beijing 100084\\ P. R. of China} \curraddr{}
\email{hqzuo@mail.tsinghua.edu.cn}

\thanks{Zuo is supported by NSFC Grant 12271280. Yau is supported by Tsinghua University Education Foundation fund (042202008).}
%\thanks{Corresponding author: Stephen S.-T. Yau.}
%\thanks{Data availability: Data sharing not applicable to this article as no datasets were generated or analysed during the current study.}

%\thanks{%Both Yau and Zuo are supported by NSFC Grant 11961141005.  Zuo is supported by NSFC Grant 12271280 and Tsinghua University Initiative Scientific Research Program.
%}

%Yau is supported by Tsinghua University start-up fund and Tsinghua University Education Foundation fund (042202008).}
%\thanks{2010 Mathematics Subject Classification:  Primary 32S25,  Secondary 58K65, 14B05.}
%\thanks{Key words: normal singularities, topological classification, weighted dual graph.}
%\date{}
 \begin{abstract}{The notion of the Yau sequence was introduced by Tomaru, as an attempt to extend Yau's elliptic sequence for (weakly) elliptic singularities to normal surface singularities of higher fundamental genera. In this paper,  we obtain the canonical cycle using the Yau cycle for certain surface singularities of degree two. Furthermore, we obtain a formula of arithmetic genera and a{\color{black} n} upper bound of  geometric genera for these singularities.  We also give some properties about the classification of weighted dual graphs of certain surface singularities of degree two.}


\end{abstract}
 \maketitle
 
 \section{Introduction}\label{sec1}

{\color{black}Let $(V, o)$ be an isolated singular point on a surface $V$, and $\pi\colon X \to V$ be a resolution of $V$. Due to the negative definiteness of the intersection matrix of $\pi^{-1}(o)$, there exists a non-zero effective divisor with support $\pi^{-1}(o)$ that has a non-positive intersection number with every exceptional curve. We define $Z$ as the unique minimal such divisor and refer to it as the \textit{fundamental cycle}. Notably, $-Z^2$ represents one of the most fundamental invariants of $(V,o)$ that remains independent of any choice in resolutions. We denote $-Z^2$ as the \textit{degree} of $(V, o)$.

In this paper, we investigate surface singularities by examining decompositions of various cycles. A primary focus is the \textit{Yau sequence}, introduced by Tomaru \cite{Tom}, which extends the first author's elliptic sequence \cite{Ya2} to singularities with larger fundamental genera. Additionally, we consider the \textit{Yau cycle} (see \cite{Ko1}), defined as the sum of all curves present in the Yau sequence (see Definition \ref{def2.14}). Notably, the Yau cycles exhibit a strong connection with canonical cycles in certain types of singularities. Konno explores applications of the Yau cycle in degree one singularities in \cite{Ko1}. In this paper, we extend Konno's results to degree two singularities and establish additional properties within a more restrictive context.}

Firstly, we discussed the relation between the canonical cycle and the Yau cycle of surface singularities of degree two. Compared to degree one situations, the degree two situations are more complicated. In  {\color{black}the} degree one situations, we know that {\color{black}the canonical cycle is obtained as a multiple of the Yau cycle}, when $Z$ is essentially irreducible. But the property doesn't hold in {\color{black} the} degree two situations, so we consider a more restrictive condition $D_m=Z_{min}$ in Theorem \ref{thm3.5} to generalize the property, where $m$ denote{\color{black}s} the length of the Yau sequence for $Z$ and $D_m$ is the smallest cycle in the Yau sequence. And we point out that the condition $D_m=Z_{min}$ holds when  $(V,o)$ are singularities of degree one.

\begin{thmletter}\label{mt1}
	Let $(V, o)$ be a normal surface singularity of degree two with $p_f(V,o)>0$, and $\pi\colon X \to V$ be a minimal resolution of $V$. Assume that fundamental cycle $Z$ is essentially irreducible and $D_m=Z_{min}$. Then $(V,o)$ is numerically Gorenstein with canonical cycle $p_f(V,o)Y$. 
\end{thmletter}

For surface singularities of degree one, Konno {\color{black} obtains} a lower bound on the arithmetic genus by using the Yau cycle $Y$ in \cite{Ko1}. For those with essentially irreducible $Z$ (see Definition \ref{def3.1}), we {\color{black} have} found that the arithmetic genus is equal to {\color{black} this} lower bound.

\begin{thmletter}\label{mt2}
Let $(V,o)$ be a normal surface singularity of degree one with $p_f(V,o)>0$, $Z$ is the fundamental cycle on the minimal resolution. Assume that $Z$ is essentially irreducible, then $p_a(V,o)=\frac{p(p-1)m}{2}+1$, where $p=p_f(V,o)$ and $m$ denotes the length of the Yau sequence for $Z$.
\end{thmletter}

We want to get a formula for the arithmetic genus of a singular point of degree 2 with essentially irreducible $Z$ like in the Theorem \ref{mt2}. For the degree two case, unlike the degree one case, if we add the additional condition that $Z$ is essentially irreducible  which still doesn't imply that $D_m=Z_{min}$. However,  the condition $D_m=Z_{min}$ is necessary when the arithmetic genus is equal to the lower bound. In order to study  singularities of degree two under the condition $D_m=Z_{min}$, In Theorem \ref{thm3.8},  we  classify the weighted dual graphs of surface singularities of the degree two when  $m>1$.  By {\color{black} classifying these graphs}, we compute  $D_m$ and $Z_{min}$ for each case in Theorem \ref{thm3.8} and obtain all cases which satisfy the conditions of $m>1$ and $D_m=Z_{min}$. From the classification of weighted dual graphs for these singularities, we obtain a formula in Theorem \ref{thm3.13} (i.e., Theorem \ref{mt3}) for  their arithmetic genera.

\begin{thmletter}\label{mt3}
Let $(V,o)$ be a normal surface singularity of degree two with $p_f(V,o)>0$, and $\pi\colon X \to V$ be a minimal resolution of $V$. Assum{\color{black} ing} that {\color{black} the} fundamental cycle $Z$ is essentially irreducible and $D_m=Z_{min}$, then 
$p_a(V,o)=[\frac{p^2}{4}]m+1$ and $p_g(V,o) \le [\frac{(p+1)^2}{4}]m$.
\end{thmletter}
 
\section{Preliminaries}\label{sec2}
\subsection{Riemann-Roch and fundamental cycle}
Let $(V, o)$ be an isolated singular point on a surface $V$, and $\pi\colon X \to V$ be a resolution of $V$. Let $\pi^{-1}(o)=\cup E_i$, $1\le i\le n$, be the decomposition of the exceptional set $\pi^{-1}(o)$ into irreducible components.

The cycles are divisors of the form $D=\sum d_iE_i$ with $d_i \in \mathbb Z$, where $E_i$ are irreducible exceptional curves. There is a natural partial ordering of the cycles: $D_1=\sum_{i} m_{i} E_{i} \le D_2=\sum_{i} n_{i} E_{i}$ if and only if $m_{i} \le n_{i}$ for all $i$. If $D_{1} \leq D_{2}$ but $D_{1} \neq D_{2}$ then we write $D_{1}<D_{2}$. We let $\textup{supp}\ D=\cup E_i$, $d_i\ne 0$, denote the support of $D$.

For a cycle $D=\sum d_iE_i$ on $\pi^{-1}(o)$, $\chi(D)$ is defined by 
$$\chi(D)=\dim H^0 ({\color{black}X},\mathcal{O}_D)-\dim H^1({\color{black}X},\mathcal{O}_D),$$
where $\mathcal{O}_D = \mathcal{O/\mathcal{O}}(-D)$. Then by Riemann-Roch theorem \cite[Proposition~IV.4, p.~75]{Se}, we have 
$$\chi (D)=-\frac12 (D^2+D\cdot K),$$
where $K$ is the canonical divisor on $X$ and $D\cdot K$ is the intersection number of $D$ and $K$. For any irreducible curve $E_i$, the adjunction formula \cite[Proposition~IV, 5, p.~75]{Se} says 
$$E_i \cdot K=-E^2_i + 2g_i +2\delta_i -2 ,$$
where $g_i$ is the genus of $E_i$ and $\delta_i$ is the degree of the conductor of $E_i$. The arithmetic genus of $D \ge 0$ is defined by $p_a(D)=1-\chi(D)$. It follows immediately from Riemann-Roch theorem that if $A$ and $B$ are cycles, then 
$$p_a(A+B)=p_a(A)+p_a(B)+A \cdot B-1.$$

\begin{definition}\label{def2.1}
Let $\pi\colon X \to V$ be a resolution, then the intersection form is negative definite on the exceptional set $\pi^{-1}(o)$. Hence, there exists a cycle $D > 0$ with support $\pi^{-1}(o)$ such that $E_i \cdot D \le 0$ for all $E_i$. We denote by $Z$ the smallest one among such cycle{\color{black}s} and call it the fundamental cycle \cite[~131-132]{Ar}.
\end{definition}

Because the fundamental cycle $Z$ is {\color{black}the} smallest, for any proper subcycle D of $Z$, there exists {\color{black}an} $E_i$ such that $E_i \cdot D >0$ and such that $E_i <Z-D$. In fact, the fundamental cycle $Z$ can be computed from the intersection as follows via a computation sequence for $Z$ in the sense of Laufer \cite[Proposition~4.1]{La2}.
\begin{align*}
Z_0=0, Z_1  &=  E_{i_1}, Z_2=Z_1+E_{i_2},\dots, Z_j=Z_{j-1}+E_{i_j},\dots,\\
Z_\ell  &=  Z_{\ell-1}+E_{i_\ell} =Z,
\end{align*}
where $E_{i_1}$ is {\color{black}an} irreducible component and $E_{i_j}\cdot Z_{j-1}>0$, $1< j\le \ell$.

Consider the computation sequence for $Z$, we have 
$$p_a(Z_j)=p_a(Z_{j-1})+p_a(E_{i_j})+Z_{j-1} \cdot E_{i_j}-1 \ge p_a(Z_{j-1})$$
for $1 \le j \le \ell$. Then $p_a(Z) \ge p_a(D)$ for any subcycle $0<D<Z$.

\begin{definition}\label{def2.2}
We {\color{black}define} the arithmetic genus of $Z$ the \textit{fundamental genus} of $(V,o)$ and denote it {\color{black}by} $p_f(V,o)$. The \textit{arithmetic genus} and \textit{geometric genus} of $(V,o)$ are respectively defined {\color{black}as} $p_a(V,o)=\max \{p_a(D)|0<D\}$ and $p_g(V,o)=\dim_C H^1(X,\mathcal{O}_X)$.
\end{definition}

\subsection{Chain-connected}we recall the notion of chain-connected {\color{black}cycles} introduced by Konno \cite{Ko2} and state {\color{black}the} fundamental properties of the chain-connected cycles.

\begin{definition}\label{def2.2}
A line bundle on a cycle is nef if it is of non-negative degree on {\color{black}all} irreducible components.
\end{definition}

\begin{definition}\label{def2.3}
A cycle $D$ is chain-connected if $\mathcal{O}_{D-C}(-C)$ is not nef for any proper subcycle $0<C<D$.
\end{definition}

According to the minimality of the fundamental cycle $Z$, the fundamental cycle is chain-connected cycle. In fact, it is the {\color{black}largest} chain-connected cycle with support $\pi^{-1}(o)$.

\begin{proposition}\label{prop2.4}
	Let D be a chain-connected cycle, then $p_a(D) \ge p_a(C)$ for any subcycle $0<C<D$.
\end{proposition}
\begin{proof}
	Since D is a chain-connected cycle, we have $\mathcal{O}_{D-A}(-A)$ is not nef for any subcycle $0<A<D$. So we can find an irreducible component $E_i<D-A$ such that $E_i \cdot A > 0$. For any subcycle $0<C<D$, like Laufer's computation sequence, we can get an increasing sequence of cycle{\color{black}s}.
Let
\begin{align*}
D_0=C, D_1  &= D_0 + E_{i_1}, D_2=D_1+E_{i_2},\dots, D_j=D_{j-1}+E_{i_j},\dots,\\
D_\ell  &=  D_{\ell-1}+E_{i_\ell} =D,
\end{align*}
where $E_{i_j}$ is {\color{black}an} irreducible component and $E_{i_j}\cdot D_{j-1}>0$, $1 \le j\le \ell$.

We have $p_a(D_j)=p_a(D_{j-1})+p_a(E_{i_j}) + E_{i_j}\cdot D_{j-1} -1 \ge p_a(D_{j-1})$ for $1 \le j \le \ell$, since $p_a(E_{i_j}) \ge 0$ and $E_{i_j}\cdot D_{j-1}>0$. We conclude that $p_a(D_\ell) \ge p_a(D_0)$, i.e., $p_a(D) \ge p_a(C)$.
\end{proof}

We remark that, when $p_a(D)=p_a(C)$, we have $p_a(D_j) = p_a(D_{j-1})$ for $1 \le j \le \ell$. It means that  $p_a(E_{i_j}) = 0$ and $E_{i_j}\cdot D_{j-1} = 1 $ for $1 \le j \le \ell$.

\begin{definition}\label{def2.5}
Let $D$ be a reducible curve. An irreducible component $E$ of $D$ is said to be a $(-m)_D-curve$ if $p_a(E)=0$ and $E \cdot (D-E) = m$.
\end{definition}

\begin{proposition}[\cite{Ko2}]\label{prop2.6}
	Given a $(-1)_D-curve$ $E$ of $D$. If D is chain-connected, then the subcycle $D'=D-E$ is {\color{black}also} chain-connected.
\end{proposition}

\begin{corollary}[\cite{Ko2}]\label{coro2.7}
	Let D be a chain-connected cycle. If $p_a(C)=p_a(D)$ for a subcycle $C<D$, then the subcycle $C$ is {\color{black}also} chain-connected.
\end{corollary}

\begin{proposition}[\cite{Ko2}]\label{prop2.8}
	Let D be a chain-connected cycle. If $\mathcal{O}_{D}(-C)$ is nef for a cycle $C$, then either $D \le C$ or $\textup{supp}\ C \cap \textup{supp}\ D= \emptyset$.
\end{proposition}

\begin{definition}\label{def2.9}
	If $D$ is chain-connected and $p_a(D) > 0$, then there uniquely exists a minimal subcycle $D_{min}$ of $D$ such that $p_a(D_{min})=p_a(D)$. We call $D_{min}$ the minimal model of $D$.
\end{definition}

\begin{theorem}[Chain-connected component decomposition, \cite{Ko2}]\label{thm2.10}
	Let $D$ be a cycle. Then there {\color{black}exists} a sequence $D_1$,$D_2$,\dots, $D_r$ of chain-connected subcycles of $D$ and a sequence $m_1$, \dots, $m_r$ of positive integers which satisfy
\begin{itemize}
\item[(1)]
$D=m_1 D_1+\dots+m_r D_r$.
\item[(2)]
For $i<j$, the cycle $-D_i$ is nef on $D_j$.
\item[(3)]
If $m_i \ge 2$, then $-D_i$ is nef on $D_i$.
\item[(4)]
For $i<j$, either $D_i > D_j$ or $\textup{supp}\ D_i \cap \textup{supp}\ D_j= \emptyset$.
\end{itemize}
\end{theorem}

\subsection{The canonical cycle}

\begin{definition}\label{def2.11}
The rational cycle $Z_{K}$ is called the canonical cycle if $Z_{K} \cdot {\color{black}E}_{i}=-K {\color{black}E}_{i}$ for all $i$, i.e.
$$
Z_{K} \cdot {\color{black}E}_{i}={\color{black}E}_{i}^{2}-2 \delta_{i}-2 g_{i}+2 \text { for all } i,
$$
where $\delta_i$ is the ``number'' of nodes and cusps on ${\color{black}E}_i$. 
\end{definition}

\begin{definition}\label{def2.12}
If the coefficients of $Z_{K}$ are integers, then the singularity is called \textit{numerical Gorenstein} {\color{black}singularity}.
\end{definition}

\subsection{Yau sequence and Yau cycle}

Assume that $p_f(V,o)=p_a(Z)>0$. By {\color{black}D}efinition {\color{black}\ref{def2.9}}, there uniquely exists a minimal subcycle $Z_{min}$ of $Z$ such that $p_a(Z_{min})=p_a(Z)$.

\begin{lemma}\label{Lemma}\label{lem2.13}
	Assume that $-Z$ is numerically trivial on $Z_{min}$. Then there uniquely exists a maximal subcycle $D<Z$ such that $\mathcal{O}_{D}(-Z)$ is numerically trivial and $p_a(D)=p_f(V,o)$. Moreover, $D$ is the fundamental cycle on its support.
\end{lemma}
\begin{proof}
Let $S=\{ 0 < D < Z | \mathcal{O}_{D}(-Z)\ is \ numerically \ trivial \ and \ p_a(D)=p_a(Z)\}$. By the assumption, $Z_{min} \in S$. Since the coefficients of subcycle D are integers, the nonempty set $S$ contains a maximal element. 

We show the maximal element of $S$ is the fundamental cycle on its support. Let $D$ be a maximal element of $S$ and $Z_1$ {\color{black}be} the fundamental cycle on $\textup{supp}\ D$. Assume that there exists an irreducible component $C \le D$ satisfying $C \cdot D > 0$. Since $\mathcal{O}_{D}(-Z)$ is numerically trivial and $C \le D$, we have $C (Z-D) < 0$ and, hence, $C \le Z-D$. Then $C+D$ is a subcycle of $Z$ and $\mathcal{O}_{C+D}(-Z)$ is numerically trivial. Furthermore, we have $p_a(Z) \ge p_a(C+D)=p_a(C)+p_a(D)+C \cdot D - 1 \ge p_a(D)=p_a(Z)$, hence $C+D \in S$. This contradicts the assumption that $D$ is maximal. Then we have $\mathcal{O}_{Z_1}(-D)$ is nef and $Z_1 \le D$. Since $D<Z$ and $p_a(D)=p_a(Z)$, the cycle D is chain-connected according to {\color{black}C}orollary  {\color{black}\ref{coro2.7}}. Furthermore, we have $D \le Z_1$ according to {\color{black}P}roposition {\color{black}\ref{prop2.8}}, since $\mathcal{O}_{D}(-Z_1)$ {\color{black}is nef} and $\textup{supp}\ D \cap \textup{supp}\ Z_1 \neq \emptyset$. Hence we get $D=Z_1$.

Assume $D_1$ and $D_2$ {\color{black}are} different maximal elements in $S$. If $\mathcal{O}_{D_2}(-D_1)$ is not nef, then there exists an irreducible component $C<D_2$ such that $C \cdot D_1 > 0$. Since $D_1$ is the fundamental cycle on its support, this shows $C \nleq D_1$ and it follows $C+D_1 \le Z$. Then we have $p_a(Z) \ge p_a(C+D_1)= p_a(C)+p_a(D_1)+C \cdot D_1 -1 \ge p_a(D_1)=p_a(Z)$ and $\mathcal{O}_{C+D_1}(-Z)$ is numerically trivial, hence $C+D \in S$. This contradicts the assumption that $D_1$ is maximal, hence we get $\mathcal{O}_{D_2}(-D_1)$ is nef. Since $D_2$ is chain-connected, either $D_2 \le D_1$ or $\textup{supp}\ D_1 \cap \textup{supp}\ D_2= \emptyset$. By the uniqueness of the minimal model $Z_{min}$, we get $Z_{min}\le D_1$ and $Z_{min} \le D_2$. It means that $D_2 \le D_1$, contradicting that $D_1$ and $D_2$ are different maximal elements in $S$. Hence there uniquely exists a maximal element in $S$.
\end{proof}

\begin{definition}[\color{black}\cite{Ko1}]\label{def2.14}
We call $D$ as in the {\color{black}L}emma {\color{black}\ref{lem2.13}} the \textit{Tyurina component} of $Z$. Since $D$ is the fundamental cycle on its support and $Z_{min}$ is also the minimal model of $D$, we can get the Tyurina component of {\color{black}$D$} when $-D$ is numerically trivial on $Z_{min}$. By the induction, we get the sequence of cycles
$$0<D_m<D_{m-1}<\dots<D_2<D_1=Z$$
such that $D_{i+1}$ is the Tyurina component of $D_i$ for $1 \le i \le m-1$ and $D_m \cdot Z_{min}<0$. We call it the \textit{Yau sequence} for $Z$ and call $Y=\sum_{i=1}^{m}D_i$ the \textit{Yau cycle}. The case $Z \cdot Z_{min}<0$ is corresponds to $m=1$.
\end{definition}

\begin{proposition}[Theorem 3.7, \cite{Ya2}]\label{prop2.15}
	If $(V,o)$ is a numerically Gorenstein elliptic singular point and $\pi \colon X \to V$ is the minimal resolution, then the Yau cycle is the canonical cycle.
\end{proposition}

The length $m$ of the Yau sequence is a numerical invariant of $(V,o)$, and gives us the arithmetic genus for singular points of fundamental genus two in \cite{Ko1}. When $p_f(V,o)>2$, length $m$ also gives us a lower bound of $p_a(V,o)$ since $p_a(Y)=\sum_{i=1}^{m}(p_a(D_i)-1)+\sum_{1 \le i<j \le m}D_i D_j +1 = m(p_f(V,o)-1)+1$.

\subsection{Classfication of weighted dual graphs}

{\color{black}Let $(V, o)$ be an isolated singular point on a surface $V$, and $\pi\colon X \to V$ be a minimal resolution of $V$. There are two beautiful results given by Artin in \cite{Ar}.
\begin{definition}\label{def2.16}
The singularity $(V, o)$ is said to be rational if $\chi (Z)=1$.
\end{definition}

If $(V, o)$ is a rational singularity, then $\pi$ is also a minimal good resolution, i.e., exceptional set with nonsingular $E_i$ and normal crossings.  Moreover, each $E_i$ is a rational curve and $E_i^2=-2$.  
\begin{theorem}[\cite{Ar}]  \label{thm2.17}   
If $(V,o)$ is a hypersurface rational singularity, then $(V,o)$ is a rational double point. Moreover the set of weighted dual graphs of hypersurface rational singularities consists of the following graphs:
\begin{flushleft}
{\Small\begin{longtable}{@{\hspace*{-25pt}}l@{\hspace*{6pt}}l@{\hspace*{12pt}}l@{\hspace*{26pt}}l@{}}
{\normalfont\upshape(1)}&$A_n, n\ge 1$&%
\scalebox{.85}{%
\begin{picture}(128,21)
\put(14,0){\line(1,0){12}} 
\put(30,0){\dashbox{2}(93,0)} 
\put(3,9){\makebox{\footnotesize$-2$}}
\put(21,9){\makebox{\footnotesize$-2$}}
\put(115,9){\makebox{\footnotesize$-2$}}
\put(11,0){\circle*{6}}
\put(29,0){\circle*{6}}
\put(123,0){\circle*{6}}
\end{picture}}
&$Z=1\ 1\dots  1$.\\[9pt]
{\upshape(2)}&$D_n, n\ge 4$&%
\scalebox{.85}{%
\begin{picture}(128,21)
\put(14,0){\line(1,0){12}} 
\put(32,0){\line(1,0){12}} 
\put(45,0){\dashbox{2}(75,0)} 
\put(3,-11){\makebox{\footnotesize$-2$}}
\put(21,-11){\makebox{\footnotesize$-2$}}
\put(39,-11){\makebox{\footnotesize$-2$}}
\put(115,-11){\makebox{\footnotesize$-2$}}
\put(11,0){\circle*{6}}
\put(29,0){\circle*{6}}
\put(29,18){\circle*{6}}
\put(31,18){\makebox{\footnotesize$-2$}}
\put(29,0){\line(0,1){18}} 
\put(47,0){\circle*{6}}
\put(29,0){\circle*{6}}
\put(123,0){\circle*{6}}
\end{picture}}
&$Z=1\ \stackrel{\textstyle 1\strut}{2\strut}\ 2\dots \ 2 \ 1.$\\[9pt]
{\upshape(3)}&$E_6$&%
\scalebox{.85}{%
\begin{picture}(128,21)
\put(14,0){\line(1,0){12}} 
\put(32,0){\line(1,0){12}} 
\put(50,0){\line(1,0){12}} 
\put(68,0){\line(1,0){12}} 
\put(3,-11){\makebox{\footnotesize$-2$}}
\put(21,-11){\makebox{\footnotesize$-2$}}
\put(39,-11){\makebox{\footnotesize$-2$}}
\put(57,-11){\makebox{\footnotesize$-2$}}
\put(75,-11){\makebox{\footnotesize$-2$}}
\put(11,0){\circle*{6}}
\put(47,18){\circle*{6}}
\put(49,18){\makebox{\footnotesize$-2$}}
\put(47,0){\line(0,1){18}} 
\put(29,0){\circle*{6}}
\put(47,0){\circle*{6}}
\put(65,0){\circle*{6}}
\put(83,0){\circle*{6}}
\end{picture}}
&$Z=1\ 2\ \stackrel{\textstyle 2\strut}{3\strut}\ 2\
1.$\\[9pt]    
{\upshape(4)}&$E_7$&
\scalebox{.85}{
\begin{picture}(128,21)
\put(14,0){\line(1,0){12}} 
\put(32,0){\line(1,0){12}} 
\put(50,0){\line(1,0){12}} 
\put(68,0){\line(1,0){12}} 
\put(86,0){\line(1,0){12}} 

\put(3,-11){\makebox{\footnotesize$-2$}}
\put(21,-11){\makebox{\footnotesize$-2$}}
\put(39,-11){\makebox{\footnotesize$-2$}}
\put(57,-11){\makebox{\footnotesize$-2$}}
\put(75,-11){\makebox{\footnotesize$-2$}}
\put(93,-11){\makebox{\footnotesize$-2$}}

\put(11,0){\circle*{6}}
\put(47,18){\circle*{6}}
\put(49,18){\makebox{\footnotesize$-2$}}
\put(47,0){\line(0,1){18}} 
\put(29,0){\circle*{6}}
\put(47,0){\circle*{6}}
\put(65,0){\circle*{6}}
\put(83,0){\circle*{6}}
\put(101,0){\circle*{6}}
\end{picture}}
&$Z=2\ 3\ \stackrel{\textstyle 2\strut}{4\strut}\ 3\ 2\
1.$\\[9pt]
{\upshape(5)}&$E_8$&
\scalebox{.85}{
\begin{picture}(128,21)
\put(14,0){\line(1,0){12}} 
\put(32,0){\line(1,0){12}} 
\put(50,0){\line(1,0){12}} 
\put(68,0){\line(1,0){12}} 
\put(3,-11){\makebox{\footnotesize$-2$}}
\put(21,-11){\makebox{\footnotesize$-2$}}
\put(39,-11){\makebox{\footnotesize$-2$}}
\put(57,-11){\makebox{\footnotesize$-2$}}
\put(75,-11){\makebox{\footnotesize$-2$}}
\put(11,0){\circle*{6}}
\put(47,18){\circle*{6}}
\put(49,18){\makebox{\footnotesize$-2$}}
\put(47,0){\line(0,1){18}} 
\put(29,0){\circle*{6}}
\put(47,0){\circle*{6}}
\put(65,0){\circle*{6}}
\put(83,0){\circle*{6}}
\put(101,0){\circle*{6}}
\put(119,0){\circle*{6}}
\put(86,0){\line(1,0){12}} 
\put(104,0){\line(1,0){12}} 
\put(93,-11){\makebox{\footnotesize$-2$}}
\put(111,-11){\makebox{\footnotesize$-2$}}
\end{picture}}
&$Z=2\ 4\ \stackrel{\textstyle 3\strut}{6\strut}\ 5\ 4\ 3\ 2.$
\end{longtable}}
\end{flushleft}
\end{theorem}

To each such weighted dual graph is associated an intersection matrix whose $(i,j)$th entry is
$E_{i} \cdot E_{j}$.

These graphs $(1)-(5)$ in Theorem \ref{thm2.17} are called ADE graphs in the literature. This theorem completely classifies the weighted dual graphs with all $E_i^2=-2$. In general, according to \cite{Ar} and \cite{Gr}, to classify the weighted dual graphs we need to classify the corresponding negative-definite matrices:
\begin{proposition}[\cite{Ar}]\label{prop2.18}
	Let $\left\{E_{i}\right\}_{i=1, \cdots, n}$ be a connected bunch of complete curves on a regular two-dimensional scheme:
	\begin{itemize}
\item[(i)] Suppose that $E_{i} \cdot E_{j}$ is negative-definite, then there exists positive cycles $Z=\sum r_{i} E_{i}$ such that $Z \cdot E_{i} \leq 0$ for all $i .$
\item[(ii)] Conversely, if there exists a positive cycle $Z=\sum r_{i} E_{i}$ such that $Z \cdot E_{i} \leq 0$ for all $i$, then $E_{i} \cdot E_{j}$ is negative semi-definite. If in addition $Z^{2}<0$, then $E_{i} \cdot E_{j}$ is negative-definite.
	\end{itemize}
	
\end{proposition} 
}

\section{Singularities of degree two} 
\subsection{The relation between the canonical cycle and the Yau cycle}
Let $(V, o)$ be an isolated singular point on a surface $V$, and $\pi\colon X \to V$ be a minimal resolution of $V$. We know the Yau cycle is the canonical cycle if $(V, o)$ is a numerically Gorenstein elliptic singular point. This property doesn't hold when $p_f(V,o) > 0$. But {\color{black}i}f the degree of $(V,o)$ is small, we can get {\color{black}a} similar property {\color{black}under} restrictive {\color{black}conditions}. If the fundamental cycle $Z$ satisfying $Z^2=-1$, the relation between the canonical cycle and the Yau cycle {\color{black}is} given by Konno in \cite{Ko1}.

\begin{definition}\label{def3.1}
Since $p_f(V,o)>0$, there exists an irreducible component $A \le Z$ such that $A$ is not {\color{black}a} $(-2)$-curve. Let $k$ be the coefficient of cycle $A$ of $Z$. We say that $Z$ is essentially irreducible if either $Z={\color{black}k}A$ or $Z-{\color{black}k}A$ consists of $(-2)$-curves.
\end{definition}

\begin{proposition}[Lemma 3.4, \cite{Ko1}]\label{prop3.2}
Let $(V, o)$ be a normal surface singularity with $p_f(V,o)>0$ and $Z^2=-1$, and $\pi\colon X \to V$ be a minimal resolution of $V$. Assume that $Z$ is essentially irreducible. Then $(V,o)$ is numerically Gorenstein with canonical cycle $(2p_f(V,o)-1)Y$.
\end{proposition}

For the case with $Z^2=-2$, we consider a similar property: $(V,o)$ is numerically Gorenstein with canonical cycle $p_f(V,o)Y$. However this property doesn't hold when only $Z$ is essentially irreducible. So we need a more restrictive situation.

\begin{lemma}\label{Lemma}\label{lem3.3}
	Let $(V, o)$ be a normal surface singularity with $p_f(V,o)>0$ and $Z^2=-2$, and $\pi\colon X \to V$ be a minimal resolution of $V$. Assume that $Z$ is essentially irreducible, then either $m=1$ or $Z-D_m$ consists of $(-2)$-curves, where $m$ denotes the length of the Yau sequence for $Z$ and $D_m$ is the smallest term in the Yau sequence for $Z$.
\end{lemma}
\begin{proof}
If $m>1$, we can get an increasing sequence of cycle{\color{black}s} like the proof of proposition 2.4:
\begin{align*}
Z_0=D_m, Z_1  &= Z_0 + E_{i_1}, Z_2=Z_1+E_{i_2},\dots, Z_j=Z_{j-1}+E_{i_j},\dots,\\
Z_\ell  &=  Z_{\ell-1}+E_{i_\ell} =Z,
\end{align*}
where $E_{i_j}$ is {\color{black}an} irreducible component and $E_{i_j}\cdot Z_{j-1}>0$, $1 \le j\le \ell$.

By the definition of Yau sequence, we have $p_a(D_m)=p_a(Z)$. Since $p_a(Z_j)=p_a(Z_{j-1})+p_a(E_{i_j}) + E_{i_j}\cdot Z_{j-1} -1 \ge p_a(Z_{j-1})$ for $1 \le j \le \ell$, we can get $p_a(E_{i_j}) = 0$ and $E_{i_j}\cdot Z_{j-1} = 1 $ for $1 \le j \le \ell$. 

Consider that $Z_{j-1} ^2 = Z_{j}^2 - 2 E_{i_j}\cdot Z_{j-1} -E_{i_j}^2=Z_{j}^2-2-E_{i_j}^2$ and $\pi$ is minimal, we have that $-2=Z_\ell^2 \le Z_{\ell-1}^2 \le \dots \le Z_1^2 \le Z_0^2 \le -1$. Assume $Z-D_m$ {\color{black}does} not consists of $(-2)$-curve, then the $E_{i_j}$ that is not $(-2)$-curve is $(-3)$-curve and only one. Since $D_m^2=-1$, we have the unique non-multiple component $A$ with $A\cdot D_m=-1$ and $\mathcal{O}_{D_m-A}(-D_m)$ is numerically trivial. Since $D_m \cdot Z_{min}<0$, $A$ is not $(-2)$-curve. Thus we can get that $\pi^{-1}(o)=\cup E_i$ satisfying all $E_i$ are $(-2)$-curves except {\color{black}for} one $(-3)$-curve $A$, and the coefficient of cycle $A$ of $Z$ is 2. Thus $p_a(Z)=(2A \cdot K-2)/2 + 1=1$ by {\color{black}the} Riemann-Roch theorem. These weighted dual graphs {\color{black}consisting of ($-2$)-curves and exactly one ($-3$)-curve} are classified in \cite{YZZ2}, there is no situation where the conditions are met. This means $Z-D_m$ consists of $(-2)$-curve{\color{black}s}.
\end{proof}

{
\color{black}
\begin{remark}
If we don't use the classification results of weighted dual graphs in \cite{YZZ2}, we can still prove Lemma \ref{lem3.3} by discussing the irreducible components with a negative intersection number with $Z$. This is the main approach to proving Theorem \ref{thm3.8} later on, however,  the process may be more cumbersome. Here, let's briefly explain the train of thought.

Continuing from the process of Lemma \ref{lem3.3}, we now need to prove that there are no singularities where the coefficient of $A$ in $Z$ is $2$ and the coefficient in $D_m$ is $1$. Since $D_m \neq Z$ and $Z$ is essentially irreducible, we know that the irreducible component with a negative intersection number with $Z$ is a ($-2$)-curve. We assert that such an irreducible component is unique. Otherwise, there exists an irreducible component $B$ such that $B \cdot Z=-1$, and the coefficient of cycle $B$ of $Z$ is $1$. It means there exists a unique irreducible component $B_1$ connected with $B$ and the coefficient of cycle $B_1$ of $Z$ is $1$. If $B_1$ is a ($-2$)-curve, then we can find the unique irreducible component $B_2$ connected with $B_1$, excluding $B$. We can continue repeating this step until $B_n$ is not longer a ($-2$)-curve. Since $Z$ is essentially irreducible, we have $B_n=A$, but the coefficient of $A$ of $Z$ is $2$, contradicting.

Then we have the unique irreducible component $C$ such that $C \cdot Z=-1$, and the coefficient of cycle $C$ of $Z$ is $2$. After removing $A$, we observe that $Z$ can be divided into several connected branches consisting of ($-2$)-curves. We will only consider the branch where $C$ is located. Its weighted dual graph is an ADE. The number of the irreducible components connected with $C$ is either $1$ or 2. By using a similar method as mentioned above, we can determine the weighted dual graph of the branch where $C$ is located is $A_n$ for some $n$ when the number is 1, and it is $D_n$ for some $n$ when the number is $2$. In these cases, we can conclude that the coefficient of $A$ of $D_m$ is $2$. The specific process can be referred to  case (2) and case (5) of Theorem \ref{thm3.8}.
\end{remark}
}

\begin{remark}\label{rm3.4}
	 It is easy to see that when $Z-D_m$ consists of $(-2)$-curves, then $D_i^2=Z^2$ for $1 \le i \le m$.
\end{remark}

\begin{theorem} [i.e. Theorem \ref{mt1}] \label{thm3.5}
	Let $(V, o)$ be a normal surface singularity with $p_f(V,o)>0$ and $Z^2=-2$, and $\pi\colon X \to V$ be a minimal resolution of $V$. Assume that $Z$ is essentially irreducible and $D_m=Z_{min}$. Then $(V,o)$ is numerically Gorenstein with canonical cycle $p_f(V,o)Y$. 
\end{theorem}
\begin{proof}
Assume $A$ is {\color{black}an} irreducible component $A\le Z$ such that A is not {\color{black}a} $(-2)$-curve, and $k$ is the coefficient of cycle $A$ of $Z$. Since $-2=Z_{min} \cdot (kA)=k Y \cdot A$, we have $K_X \cdot A=\frac{1}{k} K_X \cdot Z= \frac{2p_f(V,o)}{k}=- p_f(V,o) Y \cdot A$. We claim that $B \cdot Y=0$ for any component $B<Z-kA$. Let $j$ be the {\color{black}largest} index such that $B \le D_j$. If $j=m$, then $\mathcal{O}_{B}(-D_i)$ is numerically trivial for $1 \le i \le m-1$ by {\color{black}the} definition of {\color{black}the} Yau sequence and $B \cdot D_m=0$ by $D_m=Z_{min}$, hence $B \cdot Y=0$. If $j<m$, then $D_{j+1}$ is the Tyurina component of $D_j$. We have $\mathcal{O}_{B}(-D_i)$ is numerically trivial for $1 \le i \le j-1$ and $\mathcal{O}_{D_k}(-(D_j-D_{j+1}))$ is numerically trivial for $j+2 \le k \le m$. Since $D_k$ is chain-connected, we have either $D_k \le D_j-D_{j+1}$ or  $\textup{supp}\ D_k \cap \textup{supp}\ (D_j-D_{j+1})= \emptyset$. However, $D_k \le D_j-D_{j+1}$ is impossible, since $p_a(D_k+D_{j+1})=p_a(D_k)+p_a(D_{j+1})-1 \ge p_a(D_{j+1})$ and $\mathcal{O}_{D_k+D_{j+1}}(-D_j)$ is numerically trivial. Thus $\textup{supp}\ D_k \cap \textup{supp}\ (D_j-D_{j+1})= \emptyset$ and $C \cdot D_k=0$ for any irreducible component $C \le D_j-D_{j+1}$. In particular, $B \cdot D_k =0$ for $j+2 \le k \le m$. In sum, we get $B \cdot Y=B \cdot D_j+B \cdot D_{j+1}$. 

By Lemma \ref{lem3.3} and Remark \ref{rm3.4}, we have that $D_j^2=-2$ and $D_j-D_{j+1}$ consists of $(-2)$-curves. Since $D_j$ is the fundamental cycle on its support, there either exist two component{\color{black}s} $C_1$ and $C_2$ of $D_j$ such that $C_1 \cdot D_j=C_2\cdot D_j=-1$ and the coefficient of cycle $C_i$ of $D_j$ is 1 for $i=1,2$, or exists {\color{black}a} component $C_3$ of $D_j$ such that $C_3 \cdot D_j=-1$ and the coefficient of cycle $C_3$ of $D_j$ is 2. If the first alternative happens, then $\mathcal{O}_{D_j-C_1-C_2}(-D_j)$ is numerically trivial and $C_1 \cdot C_2=0$. This implies that $D_j-C_1-C_2$ is the Tyurina component of $Z$, and $B \in \{ C_1,C_2\}$. We have $B \cdot D_j = -1$ and $B \cdot D_{j+1}=B \cdot  (D_j-C_1-C_2)=B\cdot D_j - B^2=1$, hence $B\cdot Y=B \cdot D_j+B \cdot D_{j+1}=0$. If the last alternative happens, we can get an increasing sequence of cycle{\color{black}s} like the proof of {\color{black}P}roposition \ref{prop2.4}:
\begin{align*}
Z_0=D_{j+1}, Z_1  &= Z_0 + E_{i_1}, Z_2=Z_1+E_{i_2},\dots, Z_j=Z_{j-1}+E_{i_j},\dots,\\
Z_\ell  &=  Z_{\ell-1}+E_{i_\ell} =D_j,
\end{align*}
where $E_{i_j}$ is {\color{black}an} irreducible component and $E_{i_j}\cdot Z_{j-1}=1$ for $1 \le j\le \ell$ by $p_a(D_j)=p_a(D_{j+1})$. Then we have $E_{i_1}=C_3$ and $C_3 \cdot Y=C_3 \cdot D_j + C_3 \cdot D_{j+1}=0$. If $B \neq C_3$, it means that $B \cdot D_j=0$. Since $B \nleq D_{j+1}$ and $\mathcal{O}_{D_{j+1}+B}(-D_j)$ is numerically trivial, we have $B \cdot D_{j+1} \ge 0$ and $B \cdot D_{j+1} \le 0$, furthermore, $B \cdot Y=B \cdot D_j+B \cdot D_{j+1}=0$.

In sum, we have shown that $(V,o)$ is numerically Gorenstein with canonical cycle $p_f(V,o)Y$.
\end{proof}

We remark that, when $Z^2=-1$, the condition that $Z$ is essentially irreducible assuming that $D_m=Z_{min}$. {\color{black}However} when $Z^2=-2$, there exists {\color{black}a} singularity $(V, o)$ satisfies that $Z$ is essentially irreducible and $D_m \neq Z_{min}$.

\begin{proposition}\label{rm3.6}
	Let $(V, o)$ be a normal surface singularity with $p_f(V,o)>0$, and $\pi\colon X \to V$ be a minimal resolution of $V$. Assume that $Z$ is essentially irreducible with $K_X \cdot A + Z^2 \ge 0$ and $D_m=Z_{min}${\color{black}, where $A$ is the unique cycle in $Z$ that is not a ($-2$)-curve.} Then we have $Z_K =(\frac{2-2p_f(V,o)}{Z^2}+1) Y$. 
\end{proposition}
{\color{black}\begin{proof}
Since this property is not relevant to the main conclusions of this paper, only a brief outline of the proof is provided here.

First, let us retain the notations used in Lemma \ref{lem3.3}. We have that  $-Z^2=Z_\ell^2 \le Z_{\ell-1}^2 \le \dots \le Z_1^2 \le Z_0^2 = D_m^2 \le -1$. Since $Z$ is essentially irreducible, assume $k$ is the coefficient of cycle $A$ of $Z$, and $k'$ is the coefficient of cycle $A$ of $D_m$. Then we have that $p_a(Z)=1+\frac12(Z^2+k A\cdot K)$ and $p_a(D_m)=1+\frac12(D_m^2+k' A\cdot K)$ by the Riemann-Roch theorem. Since $Z^2 + A \cdot K\ge 0$, we have $k'=k$. It means that $Z-D_m$ consists of ($-2$)-curves and $D_i^2=Z^2$ for $1\le i \le m$.

To prove $Z_K =(\frac{2-2p_f(V,o)}{Z^2}+1) Y$, we need to demonstrate $(\frac{2-2p_f(V,o)}{Z^2}+1) Y \cdot A=-K_X \cdot A$ and $Y \cdot B=0$ for any ($-2$)-curve $B$. The former can be derived from that $Z^2=Z_{min} \cdot (kA)=k Y \cdot A$ and $p_a(Z)=1+\frac12(Z^2+k A\cdot K)$. Similar to the proof process of Theorem \ref{thm3.5}, the latter can be converted to $B \cdot (D_j+D_{j+1})=0$, where $j$ is the largest index such that $B \le D_j$. We can assume  $B \cdot D_j =-1$, since the case when $B \cdot D_j=0$ is similar to the last part of the proof of Theorem \ref{thm3.5}.  

After removing $A$, we observe that $D_j$ can be divided into several connected branches consisting of ($-2$)-curves; we will only consider the branch where $B$ is located, and claim that $B$ is the only irreducible component in this branch satisfying $B \cdot D_j < 0$. Otherwise, if we can find a path connecting two irreducible components with a negative intersection number with $D_j$, subtracting $1$ from the coefficients of the curves in this path of $D_j$ would result in a new curve with an arithmetic genus larger than that of $D_j$, contradicts the chain-connectivity of $D_j$. Additionally, as shown in the proof of Proposition \ref{prop2.4}, we can transform $D_j$ into $D_{j+1}$ by progressively removing irreducible components. We further claim that during this process, $B$ is the last irreducible component to be removed in this branch. This implies that $B \cdot D_{j+1}=1$. Thus, we have proven this remark.
\end{proof}
}

\subsection{Classification of weighted dual graphs of the singularities of degree two with the fundamental cycle $Z$ being essentially irreducible.}
Let $(V, o)$ be a normal surface singularity with $p_f(V,o)>0$ and $Z^2=-2$, and $\pi\colon X \to V$ be a minimal resolution of $V$. Assume that $Z$ is essentially irreducible and {\color{black}the} irreducible componet $A \le Z$ is not {\color{black}a} $(-2)$-curve. Let $k$ be the coefficient of cycle $A$ of $Z$ and $m$ {\color{black}be} the length of the Yau sequence for $Z$. If $m=1$, the condition $D_m=Z_{min}$ means $\mathcal{O}_{Z-kA}(-Z)$ is numerically {\color{black}trivial}. In order to study the relation between $Z$ and $D_m$, we {\color{black}aim} to classif{\color{black}y} {\color{black}the} weighted dual graphs of the singularities of degree two with the fundamental cycle $Z$ being essentially irreducible. We use $\Gamma$ to denote the weighted dual tree graph of the exceptional set $\pi^{-1}(o)$. After removing the point corresponding to the $A$, the remaining connected graphs are denoted as $\Gam_1,\cdots,\Gam_n$.

\begin{lemma} \label{lem3.7}
	With the notations as above, we have $\Gam_i$ must be ADE for any $1\le i\le n$. Let $Z|_{\Gam_i}$ be the fundamental cycle $Z$ restricted to each $\Gam_i$. If $m>1$, we have $\{ i |(Z|_{\Gam_i} \cdot Z) <0 \} =\{ i | \textup{supp} \ \Gam_i \cap \textup{supp} \ (Z-D_m) \neq \emptyset \}$, where $m$ denotes the length of the Yau sequence for $Z$ and $D_m$ denotes the smallest component of the Yau sequence.
\end{lemma}
\begin{proof}
Firstly, we have $Z|_{\Gam_i}\cdot E_j \leq 0$, $\forall E_j \in \textup{supp} \ \Gam_i$. Assume cycle $E_{j_0} \in \textup{supp} \ \Gam_i$ connected with $A$, then $Z|_{\Gam_i}\cdot E_{j_0}<0$. By Proposition \ref{prop2.18}, we conclude that $\Gam_i$ is negative-definite. Hence $\Gam_i$ must be $ADE$. 

According to Lemma \ref{lem3.3}, we have an increasing sequence of cycle{\color{black}s} when $m>1$:
\begin{align*}
Z_0=D_m, Z_1  &= Z_0 + E_{i_1}, Z_2=Z_1+E_{i_2},\dots, Z_j=Z_{j-1}+E_{i_j},\dots,\\
Z_\ell  &=  Z_{\ell-1}+E_{i_\ell} =Z,
\end{align*}
where $E_{i_j}$ is {\color{black}a} $(-2)$-curve and $E_{i_j}\cdot Z_{j-1}=1$, $1 \le j\le \ell$.

For any $1 \le i \le n$, if $Z|_{\Gam_i}\cdot Z=0$, we prove that $\mathcal{O}_{Z|_{\Gam_i}}(-Z_j)$ is numerically {\color{black}trivial} for each $j$ by induction on $j$. Firstly, since $Z|_{\Gam_i}\cdot Z=0$, we have $\mathcal{O}_{Z|_{\Gam_i}}(-Z_\ell)$ is numerically {\color{black}trivial}. Assume $\mathcal{O}_{Z|_{\Gam_i}}(-Z_j)$ is numerically {\color{black}trivial} for any $1 \le j \le \ell$, notice that $E_{i_j} \cdot Z_j=E_{i_j} \cdot Z_{j-1} + E_{i_j}^2=-1$, we have $E_{i_j} \notin \textup{supp} \ \Gam_i$ and $E_{i_j} \neq A$. It means that $E_{i_j}$ is disjoint from $Z|_{\Gam_i}$, so that $-Z_{j-1}=-(Z_j-E_{i_j})$ is numerically {\color{black}trivial} $Z|_{\Gam_i}$. By induction, we have $\mathcal{O}_{Z|_{\Gam_i}}(-Z_j)$ is numerically {\color{black}trivial} for each $j$, hence $\textup{supp} \ \Gam_i \cap \textup{supp} \ (Z-D_m) = \emptyset$.

Since $\mathcal{O}_{D_m}(-Z)$ is numerically {\color{black}trivial} when $m>1$, it is obvious that $\textup{supp} \ \Gam_i \cap \textup{supp} \ (Z-D_m) \neq \emptyset$ when $Z|_{\Gam_i} \cdot Z <0$. Hence, we have $\{ i |(Z|_{\Gam_i} \cdot Z) <0 \} =\{ i | \textup{supp} \ \Gam_i \cap \textup{supp} \ (Z-D_m) \neq \emptyset \}$.
\end{proof}

Let $S=\{ i |(Z|_{\Gam_i} \cdot Z) <0 \}$, by the proof of the Lemma \ref{lem3.7}, in order to study the relation between $Z$ and $D_m$, we only need to observe that the subgraph $\Gam'$ of $\Gam$ by removing all the $\Gam_i$, $i \notin S$. In a dual graph, the $*$ represents the point corresponding to cycle $A$.  We still call it the cycle $A$ later. The others are the points corresponding to the {\color{black}(}$-2${\color{black})-}curves, denoted by $\begin{picture}(12,3)
	\put(6,3){\circle*{6}}
\end{picture}$, we call them {\color{black}(}$-2${\color{black})-}{\color{black}curves} later.

\begin{theorem} \label{thm3.8}
With the notations as above, when $m>1$ case, $\Gam'$ and $Z|_{\Gam'}$ must be one of the following(the underlined number represents the $A$):
\begin{itemize}
\item[(1)]
$A_{m'}+A+A_{n'}$.\\
\begin{picture}(145,18)
	\put(20,-18){\makebox{\footnotesize$m'$ points}} 
	\put(16,-4){\makebox{\footnotesize$\underbrace{\hspace{40pt}}$}} 
	\put(98,-18){\makebox{\footnotesize$n'$ points}} 
	\put(94,-4){\makebox{\footnotesize$\underbrace{\hspace{40pt}}$}} 
	\put(11,3){\dashbox{2}(48,0)}
	\put(89,3){\dashbox{2}(48,0)}
	%\put(11,3){\dashbox{1}(0,24)} 
	%\put(-7,0){$*$}
	\put(71,0){$*$}
	\put(11,3){\circle*{6}}
	%\put(11,24){\circle*{6}}
	\put(59,3){\circle*{6}}
	%\put(-4,3){\line(1,0){12}} 
	\put(61,3){\line(1,0){12}} 
	\put(74,3){\line(1,0){12}}
	\put(89,3){\circle*{6}}
	\put(137,3){\circle*{6}}
\end{picture}
$Z|_{\Gam'} =1 \ ... \ \underline{1} \ ... \ 1$.
\\
\\
\item[(2)]
$A+A_{n'}$: $n'\ge 3$.\\
\begin{picture}(83,18)
	\put(20,-18){\makebox{\footnotesize$n'$ points}} 
	\put(16,-4){\makebox{\footnotesize$\underbrace{\hspace{40pt}}$}} 
	\put(11,3){\dashbox{2}(48,0)}
	%\put(11,3){\dashbox{1}(0,24)} 
	%\put(-7,0){$*$}
	\put(71,0){$*$}
	\put(11,3){\circle*{6}}
	%\put(11,24){\circle*{6}}
	\put(59,3){\circle*{6}}
	%\put(-4,3){\line(1,0){12}} 
	\put(61,3){\line(1,0){12}} 
\end{picture}
$Z|_{\Gam'} =1 \ 2 \ 2 \ ...\ 2 \ \underline{2}$,or 
\begin{picture}(83,18)
	\put(20,-18){\makebox{\footnotesize$n'$ points}} 
	\put(16,-4){\makebox{\footnotesize$\underbrace{\hspace{40pt}}$}} 
	\put(11,3){\dashbox{2}(48,0)}
	%\put(11,3){\dashbox{1}(0,24)} 
	%\put(-7,0){$*$}
	\put(71,0){$*$}
	\put(11,3){\circle*{6}}
	%\put(11,24){\circle*{6}}
	\put(59,3){\circle*{6}}
	%\put(-4,3){\line(1,0){12}} 
	\put(61,5){\line(1,0){12}} 
	\put(61,1){\line(1,0){12}} 
\end{picture}
$Z|_{\Gam'} =1 \ 2 \ 2 \ ...\ 2 \ \underline{1}$.
\\
\\
\item[(3)]
$A+(1-A_{n'})$: $n'\ge 3$.\\
\begin{picture}(83,18)
	%\put(20,-18){\makebox{\footnotesize$n-1$ points}} 
	%\put(16,-4){\makebox{\footnotesize$\underbrace{\hspace{40pt}}$}} 
	\put(11,3){\dashbox{2}(48,0)}
	%\put(11,3){\dashbox{1}(0,24)} 
	%\put(-7,0){$*$}
	\put(71,0){$*$}
	\put(11,3){\circle*{6}}
	%\put(11,24){\circle*{6}}
	\put(59,3){\circle*{6}}
	%\put(-4,3){\line(1,0){12}} 
	\put(61,3){\line(1,0){12}} 
	\put(59,5){\line(0,1){12}} 
	\put(59,19){\circle*{6}}
\end{picture}
$Z|_{\Gam'} = 1 \ 2 \ ...\stackrel {\textstyle {1}\strut}{2\strut} \ \underline{1}$.
\\
\item[(4)]
$A+(k'-D_{n'})$: k' is an even number and $0\le k' \le n-3$.\\
\begin{picture}(131,18)
	\put(20,-18){\makebox{\footnotesize$k'+1$ points}} 
	\put(16,-4){\makebox{\footnotesize$\underbrace{\hspace{40pt}}$}} 
	\put(11,3){\dashbox{2}(48,0)}
	\put(59,3){\dashbox{2}(48,0)}
	%\put(11,3){\dashbox{1}(0,24)} 
	%\put(-7,0){$*$}
	\put(107,3){\circle*{6}}
	\put(123,3){\circle*{6}}
	\put(107,19){\circle*{6}}
	\put(11,3){\circle*{6}}
	%\put(11,24){\circle*{6}}
	\put(59,3){\circle*{6}}
	%\put(-4,3){\line(1,0){12}} 
	\put(109,3){\line(1,0){12}} 
	\put(107,5){\line(0,1){12}} 
	\put(59,5){\line(0,1){12}} 
	\put(56,16){$*$}
\end{picture}
$Z|_{\Gam'} = 2 \ 3 \ ...\stackrel {\textstyle \underline{1}\strut}{k+2\strut} \ ... \stackrel {\textstyle {\frac{k'+2}{2}}\strut}{k'+2\strut} \frac{k'+2}{2}$.
\\
\\
\\
In particular, when $n'$ is an odd number and $k'=n'-3$, we have:\\
\begin{picture}(83,18)
	%\put(20,-18){\makebox{\footnotesize$k+1$ points}} 
	%\put(16,-4){\makebox{\footnotesize$\underbrace{\hspace{40pt}}$}} 
	\put(11,3){\dashbox{2}(48,0)}
	%\put(11,3){\dashbox{1}(0,24)} 
	%\put(-7,0){$*$}
	\put(75,3){\circle*{6}}
	\put(59,-13){\circle*{6}}
	\put(11,3){\circle*{6}}
	%\put(11,24){\circle*{6}}
	\put(59,3){\circle*{6}}
	%\put(-4,3){\line(1,0){12}} 
	\put(61,3){\line(1,0){12}} 
	\put(59,1){\line(0,-1){12}} 
	\put(59,5){\line(0,1){12}} 
	\put(56,16){$*$}
\end{picture}
$Z|_{\Gam'} = 2 \ 3 \ ...\underset{\textstyle \frac{n'-1}{2}} {\stackrel {\textstyle \underline{1}\strut}{n'-1\strut}}  \frac{n'-1}{2}$.\\
In additional, when $k'=0$, we have:\\
\begin{picture}(99,18)
	%\put(20,-18){\makebox{\footnotesize$k+1$ points}} 
	%\put(16,-4){\makebox{\footnotesize$\underbrace{\hspace{40pt}}$}} 
	\put(27,3){\dashbox{2}(48,0)}
	%\put(11,3){\dashbox{1}(0,24)} 
	%\put(-7,0){$*$}
	\put(91,3){\circle*{6}}
	\put(75,19){\circle*{6}}
	\put(27,3){\circle*{6}}
	%\put(11,24){\circle*{6}}
	\put(75,3){\circle*{6}}
	%\put(-4,3){\line(1,0){12}} 
	\put(11,3){\line(1,0){12}} 
	\put(77,3){\line(1,0){12}} 
	\put(75,5){\line(0,1){12}} 
	\put(7,0){$*$}
\end{picture}
$Z|_{\Gam'} =\underline{1} \ 2 \ 2 \ ... \stackrel {\textstyle 1 \strut}{2\strut}  1$.
\\
\item[(5)]
$A+(D_{n'}')$: $n'$ is an odd number.\\
\begin{picture}(99,18)
	%\put(20,-18){\makebox{\footnotesize$k+1$ points}} 
	%\put(16,-4){\makebox{\footnotesize$\underbrace{\hspace{40pt}}$}} 
	\put(11,3){\dashbox{2}(48,0)}
	%\put(11,3){\dashbox{1}(0,24)} 
	%\put(-7,0){$*$}
	\put(75,3){\circle*{6}}
	\put(59,19){\circle*{6}}
	\put(11,3){\circle*{6}}
	%\put(11,24){\circle*{6}}
	\put(59,3){\circle*{6}}
	%\put(-4,3){\line(1,0){12}} 
	\put(61,3){\line(1,0){12}} 
	\put(59,5){\line(0,1){12}} 
	\put(77,3){\line(1,0){12}} 
	\put(87,0){$*$}
\end{picture}
$Z|_{\Gam'} = 2 \ 3 \ ...\stackrel{\textstyle \frac{n'-1}{2} \strut}{n'-1\strut} \ \frac{n'+1}{2} \ \underline{2}$,\\
or \begin{picture}(99,18)
	%\put(20,-18){\makebox{\footnotesize$k+1$ points}} 
	%\put(16,-4){\makebox{\footnotesize$\underbrace{\hspace{40pt}}$}} 
	\put(11,3){\dashbox{2}(48,0)}
	%\put(11,3){\dashbox{1}(0,24)} 
	%\put(-7,0){$*$}
	\put(75,3){\circle*{6}}
	\put(59,19){\circle*{6}}
	\put(11,3){\circle*{6}}
	%\put(11,24){\circle*{6}}
	\put(59,3){\circle*{6}}
	%\put(-4,3){\line(1,0){12}} 
	\put(61,3){\line(1,0){12}} 
	\put(59,5){\line(0,1){12}} 
	\put(77,5){\line(1,0){12}} 
	\put(77,1){\line(1,0){12}} 
	\put(87,0){$*$}
\end{picture}
$Z|_{\Gam'} = 2 \ 3 \ ...\stackrel{\textstyle \frac{n'-1}{2} \strut}{n'-1\strut} \ \frac{n'+1}{2} \ \underline{1}$.\\
\item[(6)]
$A+(D_{n'}'')$: $n'$ is an even number.\\
\begin{picture}(83,18)
	%\put(20,-18){\makebox{\footnotesize$k+1$ points}} 
	%\put(16,-4){\makebox{\footnotesize$\underbrace{\hspace{40pt}}$}} 
	\put(11,3){\dashbox{2}(48,0)}
	%\put(11,3){\dashbox{1}(0,24)} 
	%\put(-7,0){$*$}
	\put(75,3){\circle*{6}}
	\put(59,19){\circle*{6}}
	\put(11,3){\circle*{6}}
	%\put(11,24){\circle*{6}}
	\put(59,3){\circle*{6}}
	%\put(-4,3){\line(1,0){12}} 
	\put(61,3){\line(1,0){12}} 
	\put(59,5){\line(0,1){12}} 
	\put(75,5){\line(0,1){12}} 
	\put(61,19){\line(1,0){12}} 
	\put(72,16){$*$}
\end{picture}
$Z|_{\Gam'} = 2 \ 3 \ ...\stackrel{\textstyle \frac{n'}{2} \strut}{n'-1\strut}\ \stackrel {\textstyle \underline{1} \strut}{\frac{n'}{2}\strut}$.\\ 
\item[(7)]
$A+E_6$.\\ 
\begin{picture}(98,18)
	\put(11,3){\circle*{6}}
	\put(11,3){\line(1,0){12}}
	\put(27,3){\circle*{6}}
	\put(27,3){\line(1,0){12}}
	\put(43,3){\circle*{6}}
	\put(43,3){\line(1,0){12}}
	\put(59,3){\circle*{6}}
	\put(59,3){\line(1,0){12}}
	\put(43,5){\line(0,1){12}}
	\put(43,19){\circle*{6}}
	\put(75,3){\circle*{6}}
	\put(75,3){\line(1,0){12}}
	\put(86,0){$*$}
\end{picture}	
$Z|_{\Gam'} = 2 \ 3 \stackrel{\textstyle 2 \strut}{4\strut} 3 \ 2 \ \underline{1}$.\\ 

\item[(8)]
$A+D_5'''$: \\
\begin{picture}(82,18)
	\put(11,3){\circle*{6}}
	\put(11,3){\line(1,0){12}}
	\put(27,3){\circle*{6}}
	\put(27,3){\line(1,0){12}}
	\put(43,3){\circle*{6}}
	\put(43,3){\line(1,0){12}}
	\put(59,3){\circle*{6}}
	\put(59,3){\line(1,0){12}}
	\put(43,5){\line(0,1){12}}
	\put(43,19){\circle*{6}}
	\put(70,0){$*$}
\end{picture}	
$Z|_{\Gam'} = 1 \ 2 \stackrel{\textstyle 2 \strut}{3\strut} 2 \ \underline{1}$.

	\end{itemize}
\end{theorem}
{\color{black}The} weighted dual graph of case (8) is the same as case (5), but $Z|_{\Gam'}$ {\color{black}differs}.
\begin{proof}
Since $m>1$, $Z \cdot A=0$ and $-2=Z^2=\overset{n}{\underset{i=1}\sum} (Z \cdot Z|_{\Gam_i}) $, hence, the number of elements in $S$ is {\color{black}either} $1$ or $2$. If $S$ has two elements, we can assume $S=\{ 1,2 \}$. Then we know that $Z \cdot Z|_{\Gam_1}=-1$ and, hence, there exists {\color{black}an} irreducible componet $C_1 \le Z|_{\Gam_1}$ such that $C_1 \cdot Z=-1$ and the coefficient of $C_1$ of $Z$ is $1$. Since $C_1$ is {\color{black}a} $(-2)$-curve, there exists {\color{black}a} unique cycle in $\Gam$ connected with $C_1$. If the cycle isn't $A$, we call it $C_2$. Then we have $C_2 \cdot Z=-1$ and the coefficient of $C_1$ of $Z$ is $1$. Except for $C_1$, there exists {\color{black}a} unique cycle in $\Gam$ connected with $C_2$. Now{\color{black}, by} the obvious induction{\color{black}, we can show} that the weighted dual graph of $\Gam_1$ is $A_{m'}$ for any $m'$ and similarly the weighted dual graph of $\Gam_2$ is $A_{n'}$ for any $n'$, and the weighted dual graph of $\Gam'$ and $Z|_{\Gam'}$ is the same as {\color{black}in} case (1).

Then we consider that $S$ has only one element, we can assume $S=\{ 1\}$. Then we know that $Z \cdot Z|_{\Gam_1}=-2$ and there exists {\color{black}a} unique cycle $C$ in $\Gam_1$ such that $C \cdot Z<0$. That is because if there {\color{black}exist} two cycles $C$ and $D$, we can find a connection way of $C$ and $D$ in $\Gam_1$. We denote the points of the connection way {\color{black}by} $C_1=C$,$C_2$,$\dots$,$C_{\color{black}k'}=D$, then $\mathcal{O}_{Z}(-(Z-\overset{k'}{\underset{i=1}\sum} C_i))$ is nef, contradicting the minimality of $Z$. Since $C {\color{black}\cdot} Z=p_a(Z)-p_a(Z-C)-1 \ge -1$, the coefficient of $C$ of $Z$ is $2$.

Refer{\color{black}ring} to \cite{Oku}, for an ADE graph, with abuse of notations, $E_i$ denotes exceptional curve, and $E= \sum E_i$ is the exceptional cycle. Let $\delta_{i}=\left(E-E_{i}\right) \cdot E_{i}$ be the number of irreducible components of $E$ connected with $E_{{\color{black}i}}$. A cycle $E_{i}$ is called an end (resp. a node), if $\delta_{i}=1$ (resp. $\delta_{i} \geq 3$). 

If $C$ isn't an end in $\Gam_1$, we claim that $\Gam_1$ is $A_{n'}$ for any $n'$. If not, then $\Gam_1$ has a node $D$, we can find a connection way of $C$ and $D$. We denote the points of the connection way {\color{black}by} $C_1=C$,$C_2$,$\dots$,$C_k'=D$(If $C=D$, denote $k'=1$). If $C \neq D$, we denote the cycles connected with $D$ in $\Gam_1$ {\color{black}by} $C_{k'-1}$,$B_1$,$B_2$, and the another cycle connected with $C$ in $\Gam_1$ {\color{black}as $B_3$}. If $C=D$, we denote the cycles connected with $D$ in $\Gam_1$ as $B_1$,$B_2$,$B_3$. Then we have $\mathcal{O}_{Z}(-(Z-2\overset{k'}{\underset{i=1}\sum} C_i-\overset{3}{\underset{i=1}\sum} B_i))$ is nef, contradicting the minimality of $Z$. Since the coefficient of $C$ of $Z$ is $2$ and $C$ isn't {\color{black}an} end in $\Gam_1$, we have there exists a cycle $D$ connected with $C$ in $\Gam_1$ such that the coefficient of $D$ of $Z$ is $1$, hence, $D$ is an end. Then either $C$ {\color{black}is} connected with $A$ or the coefficient of the another cycle $C_2$ that connected with $C$ in $\Gam_1$ of $Z$ is $2$. If the first alternative happens, {\color{black}there are two situations. One is that} $C$ {\color{black}is} connected with another cycle $C_2$ and $C_2$ is an end when $kA \cdot C=1$, {\color{black}the other is that} $C$ is an end in $\Gam_1$ when $kA \cdot C=2$({\color{black}which} contradicts the assumption that $C$ isn't an end). If the last alternative happens, for $C_2$, either $C_2$ {\color{black}is} connected with $A$ {\color{black},} or the coefficient of the another cycle $C_3$ that {\color{black}is} connected with $C_2$ in $\Gam_1$ of $Z$ is $2$. By induction, we can get case (2) and case (3).

If $C$ is an end, we have $\Gam_1$ isn't $A_{n'}$ for any $n'$ by $m>1$. Then $\Gam_1$ has a node {\color{black}called} $D$, we can find a connection way of $C$ and $D$, and denote the points of the connection way {\color{black}by} $C_1=C$,$C_2$,$\dots$,$C_k'=D$, and {\color{black}denote} $B_1$,$B_2$ as the other cycles connected with $D$. Assume $c_i$(resp. $b_i$) is the coefficient of $C_i$(resp. $B_i$) of $Z$. By induction, we have $c_{j}=j+1-\overset{j-1}{\underset{i=1}\sum}((j-i)kA\cdot C_i)$ for $2\le j \le k'$. Since $b_1,b_2\ge \frac{c_k'}{2}$, we have $0 \le c_{k'}-c_{k'-1}-kA\cdot C_{k'}=1-\overset{k'}{\underset{i=1}\sum}(kA\cdot C_i)$. If $A \cdot \overset{k'}{\underset{i=1}\sum}C_i>0$, we have $B_1$,$B_2$ are {\color{black}the} ends and $k=1$, then we get case (4).

If $A \cdot \overset{k'}{\underset{i=1}\sum}C_i=0$, we have $c_j=j+1$ for $1\le j \le k'$. If $k'$ is an odd number, then we can assume $b_1=\frac{k'+1}{2}$ and $b_2=\frac{k'+3}{2}$, hence, $B_1$ is an end and $kA \cdot B_2 \le 2$. We can get case (5) when $kA \cdot B_2=2$, and get case (7) when $kA \cdot B_2=0$. If $kA \cdot B_2=1$, there exists a cycle $B_3$ connected with $B_2$, and the coefficient of $B_3$ of $Z$ is 1. But $B_3 \cdot Z \ge -2+\frac{k'+3}{2} >0$, it's impossible.

If $k'$ is an even number, then we can assume $b_1=b_2=\frac{k'+2}{2}$. We can get case (6) when $kA \cdot B_1>0$ and  $kA \cdot B_2>0$, and get case (8) when $kA \cdot B_1=0$ or  $kA \cdot B_2=0$. If $kA \cdot B_1=0$ and $kA \cdot B_2=0$, we can prove that $A \cdot Z|_{\Gam_1}=0$, it's impossible.
\end{proof}

\begin{proposition} \label{prop3.9}
	With the notations as above, when $m>1$ and $D_m=Z_{min}$ case, $\Gam'$ must be one of the following:
\begin{itemize}
\item[(1)]
$A_{n'}+A+A_{n'}$.
\item[(2)]
$A+A_{n'}$: $n'$ is an odd number and $n'\ge 3$.
\item[(3)]
$A+(1-A_{n'})$: $n'$ is an odd number and $n'\ge 3$.
\item[(4)]
$A+(k'-D_{n'})$: k' is an even number and $0\le k' \le n-3$.
\item[(5)]
$A+(D_{n'}')$: $n'$ is an odd number.
\item[(6)]
$A+(D_{n'}'')$: $n'$ is an even number.
\item[(7)]
$A+E_6$.
\end{itemize}
remark: The weighted dual graphs $\Gam'$ of case{\color{black}s} (1)-(7) in the Proposition \ref{prop3.9} {\color{black}are} the same as {\color{black}those of }case{\color{black}s} (1)-(7) {\color{black}mentioned} in the Theorem \ref{thm3.8}.
\end{proposition}
\begin{proof}
We only need to compare $D_m|_{\Gam'}$ with $Z_{min}|_{\Gam'}$, since $\textup{supp} \ (Z-D_m) \subset \textup{supp} \ (Z-Z_{min}) \subset \textup{supp} \ \Gam'$. We can compute the $D_m|_{\Gam'}$ and $Z_{min}|_{\Gam'}$ for each case in the Theorem \ref{thm3.8}. Computations {\color{black}for} case{\color{black}s} (1)-(8) are simple, {\color{black}so} let us take the computation of the case (1) as an example here. 

In case (1), there {\color{black}exist} two different irreducible components $B_1$ and $C_1$ in $D_1=Z$ such that $B_1 \cdot D_1=-1${\color{black},} $C_1 \cdot D_1=-1$, and $B_1 \cdot C_1 = 0$. In fact, the points corresponding to the $B_1$ and $C_1$ are the ends of $A_{m'}$ and $A_{n'}$ that {\color{black}are} not connect{\color{black}ed} with the point corresponding to the $A$. Then we have $D_2=D_1-B_1-C_1$. We can assume $m' \ge n'${\color{black}. B}y induction, we know that $m=n'+1$ and the weighted dual {\color{black}graph} of $D_{m}|_{\Gam'}$ {\color{black}consists} of $A_{m'-n'}$ and the point corresponding to the $A$. If $m'-n' \neq 0$, since the cycle $C$ corresponding to the end of $A_{m'-n'}$ that {\color{black}is} not connect{\color{black}ed} with the point corresponding to the $A${\color{black}, and it} satisfies $C \cdot D_m=-1$, we have $p_a(D_m)=p_a(D_m-C)$. By induction, we have $Z_{min}|_{\Gam'}=A$.

Now we give the $D_m|_{\Gam'}$ and $Z_{min}|_{\Gam'}$ for each case in the Theorem \ref{thm3.8}.
\begin{itemize}
\item[(1)] 
$A_{m'}+A+A_{n'}$: assume $m' \ge n'$.\\
$D_m|_{\Gam'}$:\begin{picture}(83,18)
	\put(20,-18){\makebox{\footnotesize$m'-n'$ points}} 
	\put(16,-4){\makebox{\footnotesize$\underbrace{\hspace{40pt}}$}} 
	\put(11,3){\dashbox{2}(48,0)}
	%\put(11,3){\dashbox{1}(0,24)} 
	%\put(-7,0){$*$}
	\put(71,0){$*$}
	\put(11,3){\circle*{6}}
	%\put(11,24){\circle*{6}}
	\put(59,3){\circle*{6}}
	%\put(-4,3){\line(1,0){12}} 
	\put(61,3){\line(1,0){12}} 
\end{picture} $1 \ ... \ \underline{1}$ .\\
\\
$Z_{min}|_{\Gam'}$:\begin{picture}(23,18)
	\put(11,0){$*$}
\end{picture} \underline{1}.\\
\item[(2)]
$A+A_{n'}$: $n'\ge 3$.\\
$D_m|_{\Gam'}$(when n' is an even number):
\begin{picture}(50,18)
	\put(11,3){\circle*{6}}
	\put(11,3){\line(1,0){12}}
	\put(27,3){\circle*{6}}
	\put(27,3){\line(1,0){12}}
	\put(38,0){$*$}
\end{picture} $1 \ 2 \ \underline{2}$
or 
\begin{picture}(50,18)
	\put(11,3){\circle*{6}}
	\put(11,3){\line(1,0){12}}
	\put(27,3){\circle*{6}}
	\put(27,5){\line(1,0){12}}
	\put(27,1){\line(1,0){12}}
	\put(38,0){$*$}
\end{picture}$1 \ 2 \ \underline{1}$.\\
$D_m|_{\Gam'}$(when n' is an odd number):
\begin{picture}(34,18)
	\put(11,3){\circle*{6}}
	\put(11,3){\line(1,0){12}}
	\put(22,0){$*$}
\end{picture} $1 \ \underline{2}$
or 
\begin{picture}(34,18)
	\put(11,3){\circle*{6}}
	\put(11,5){\line(1,0){12}}
	\put(11,1){\line(1,0){12}}
	\put(22,0){$*$}
\end{picture}$1 \ \underline{1}$.\\
$Z_{min}|_{\Gam'}$:
\begin{picture}(34,18)
	\put(11,3){\circle*{6}}
	\put(11,3){\line(1,0){12}}
	\put(22,0){$*$}
\end{picture} $1 \ \underline{2}$
or 
\begin{picture}(34,18)
	\put(11,3){\circle*{6}}
	\put(11,5){\line(1,0){12}}
	\put(11,1){\line(1,0){12}}
	\put(22,0){$*$}
\end{picture}$1 \ \underline{1}$.\\
\item[(3)]
$A+(1-A_{n'})$: $n'\ge 3$.\\
$D_m|_{\Gam'}$(when n' is an even number):
\begin{picture}(34,18)
	\put(11,3){\circle*{6}}
	\put(11,3){\line(1,0){12}}
	\put(11,5){\line(0,1){12}}
	\put(11,19){\circle*{6}}
	\put(22,0){$*$}
\end{picture} $ \stackrel{\textstyle 1 \strut}{1\strut} \underline{1}$.\\
$D_m|_{\Gam'}$(when n' is an odd number):
\begin{picture}(23,18)
	\put(11,0){$*$}
\end{picture} \underline{1}.\\
$Z_{min}|_{\Gam'}$:
\begin{picture}(23,18)
	\put(11,0){$*$}
\end{picture} \underline{1}.\\
\item[(4)]
$A+(k'-D_{n'})$: k' is an even number and $0\le k' \le n-3$.\\
$D_m|_{\Gam'}$:\begin{picture}(131,18)
	\put(20,-18){\makebox{\footnotesize$k'$ points}} 
	\put(16,-4){\makebox{\footnotesize$\underbrace{\hspace{40pt}}$}} 
	\put(11,3){\dashbox{2}(48,0)}
	\put(59,3){\dashbox{2}(48,0)}
	%\put(11,3){\dashbox{1}(0,24)} 
	%\put(-7,0){$*$}
	\put(107,3){\circle*{6}}
	\put(123,3){\circle*{6}}
	\put(107,19){\circle*{6}}
	\put(11,3){\circle*{6}}
	%\put(11,24){\circle*{6}}
	\put(59,3){\circle*{6}}
	%\put(-4,3){\line(1,0){12}} 
	\put(109,3){\line(1,0){12}} 
	\put(107,5){\line(0,1){12}} 
	\put(59,5){\line(0,1){12}} 
	\put(56,16){$*$}
\end{picture}
$1 \ 2 \ ...\stackrel {\textstyle \underline{1}\strut}{k\strut} \ ... \stackrel {\textstyle {\frac{k'}{2}}\strut}{k'\strut} \frac{k'}{2}$.
\\
\\
\\
$Z_{min}|_{\Gam'}$:\begin{picture}(131,18)
	\put(20,-18){\makebox{\footnotesize$k'$ points}} 
	\put(16,-4){\makebox{\footnotesize$\underbrace{\hspace{40pt}}$}} 
	\put(11,3){\dashbox{2}(48,0)}
	\put(59,3){\dashbox{2}(48,0)}
	%\put(11,3){\dashbox{1}(0,24)} 
	%\put(-7,0){$*$}
	\put(107,3){\circle*{6}}
	\put(123,3){\circle*{6}}
	\put(107,19){\circle*{6}}
	\put(11,3){\circle*{6}}
	%\put(11,24){\circle*{6}}
	\put(59,3){\circle*{6}}
	%\put(-4,3){\line(1,0){12}} 
	\put(109,3){\line(1,0){12}} 
	\put(107,5){\line(0,1){12}} 
	\put(59,5){\line(0,1){12}} 
	\put(56,16){$*$}
\end{picture}
$1 \ 2 \ ...\stackrel {\textstyle \underline{1}\strut}{k\strut} \ ... \stackrel {\textstyle {\frac{k'}{2}}\strut}{k'\strut} \frac{k'}{2}$.
\\
\\
\\
\item[(5)]
$A+(D_{n'}')$: $n'$ is an odd number.\\
$D_m|_{\Gam'}$:\begin{picture}(99,18)
	\put(20,-18){\makebox{\footnotesize$n'-3$ points}} 
	\put(16,-4){\makebox{\footnotesize$\underbrace{\hspace{40pt}}$}} 
	\put(11,3){\dashbox{2}(48,0)}
	%\put(11,3){\dashbox{1}(0,24)} 
	%\put(-7,0){$*$}
	\put(75,3){\circle*{6}}
	\put(59,19){\circle*{6}}
	\put(11,3){\circle*{6}}
	%\put(11,24){\circle*{6}}
	\put(59,3){\circle*{6}}
	%\put(-4,3){\line(1,0){12}} 
	\put(61,3){\line(1,0){12}} 
	\put(59,5){\line(0,1){12}} 
	\put(77,3){\line(1,0){12}} 
	\put(87,0){$*$}
\end{picture}
$1 \ ...\stackrel{\textstyle \frac{n'-3}{2} \strut}{n'-3\strut} \ \frac{n'-1}{2} \ \underline{2}$\\
\\
\\
or \begin{picture}(99,18)
	\put(20,-18){\makebox{\footnotesize$n'-3$ points}} 
	\put(16,-4){\makebox{\footnotesize$\underbrace{\hspace{40pt}}$}} 
	\put(11,3){\dashbox{2}(48,0)}
	%\put(11,3){\dashbox{1}(0,24)} 
	%\put(-7,0){$*$}
	\put(75,3){\circle*{6}}
	\put(59,19){\circle*{6}}
	\put(11,3){\circle*{6}}
	%\put(11,24){\circle*{6}}
	\put(59,3){\circle*{6}}
	%\put(-4,3){\line(1,0){12}} 
	\put(61,3){\line(1,0){12}} 
	\put(59,5){\line(0,1){12}} 
	\put(77,5){\line(1,0){12}} 
	\put(77,1){\line(1,0){12}} 
	\put(87,0){$*$}
\end{picture}
$1 \ ...\stackrel{\textstyle \frac{n'-3}{2} \strut}{n'-3\strut} \ \frac{n'-1}{2} \ \underline{1}$.\\
\\
\\
$Z_{min}|_{\Gam'}$:\begin{picture}(99,18)
	\put(20,-18){\makebox{\footnotesize$n'-3$ points}} 
	\put(16,-4){\makebox{\footnotesize$\underbrace{\hspace{40pt}}$}} 
	\put(11,3){\dashbox{2}(48,0)}
	%\put(11,3){\dashbox{1}(0,24)} 
	%\put(-7,0){$*$}
	\put(75,3){\circle*{6}}
	\put(59,19){\circle*{6}}
	\put(11,3){\circle*{6}}
	%\put(11,24){\circle*{6}}
	\put(59,3){\circle*{6}}
	%\put(-4,3){\line(1,0){12}} 
	\put(61,3){\line(1,0){12}} 
	\put(59,5){\line(0,1){12}} 
	\put(77,3){\line(1,0){12}} 
	\put(87,0){$*$}
\end{picture}
$1 \ ...\stackrel{\textstyle \frac{n'-3}{2} \strut}{n'-3\strut} \ \frac{n'-1}{2} \ \underline{2}$\\
\\
\\
or \begin{picture}(99,18)
	\put(20,-18){\makebox{\footnotesize$n'-3$ points}} 
	\put(16,-4){\makebox{\footnotesize$\underbrace{\hspace{40pt}}$}} 
	\put(11,3){\dashbox{2}(48,0)}
	%\put(11,3){\dashbox{1}(0,24)} 
	%\put(-7,0){$*$}
	\put(75,3){\circle*{6}}
	\put(59,19){\circle*{6}}
	\put(11,3){\circle*{6}}
	%\put(11,24){\circle*{6}}
	\put(59,3){\circle*{6}}
	%\put(-4,3){\line(1,0){12}} 
	\put(61,3){\line(1,0){12}} 
	\put(59,5){\line(0,1){12}} 
	\put(77,5){\line(1,0){12}} 
	\put(77,1){\line(1,0){12}} 
	\put(87,0){$*$}
\end{picture}
$1 \ ...\stackrel{\textstyle \frac{n'-3}{2} \strut}{n'-3\strut} \ \frac{n'-1}{2} \ \underline{1}$.\\
\\
\item[(6)]
$A+(D_{n'}'')$: $n'$ is an even number.\\
$D_m|_{\Gam'}$:\begin{picture}(83,18)
	\put(20,-18){\makebox{\footnotesize$n'-3$ points}} 
	\put(16,-4){\makebox{\footnotesize$\underbrace{\hspace{40pt}}$}} 
	\put(11,3){\dashbox{2}(48,0)}
	%\put(11,3){\dashbox{1}(0,24)} 
	%\put(-7,0){$*$}
	\put(75,3){\circle*{6}}
	\put(59,19){\circle*{6}}
	\put(11,3){\circle*{6}}
	%\put(11,24){\circle*{6}}
	\put(59,3){\circle*{6}}
	%\put(-4,3){\line(1,0){12}} 
	\put(61,3){\line(1,0){12}} 
	\put(59,5){\line(0,1){12}} 
	\put(75,5){\line(0,1){12}} 
	\put(61,19){\line(1,0){12}} 
	\put(72,16){$*$}
\end{picture}
$1 \ ...\stackrel{\textstyle \frac{n'-2}{2} \strut}{n'-3\strut}\ \stackrel {\textstyle \underline{1} \strut}{\frac{n'-2}{2}\strut}$.
\\ 
\\
\\
$Z_{min}|_{\Gam'}$:\begin{picture}(83,18)
	\put(20,-18){\makebox{\footnotesize$n'-3$ points}} 
	\put(16,-4){\makebox{\footnotesize$\underbrace{\hspace{40pt}}$}} 
	\put(11,3){\dashbox{2}(48,0)}
	%\put(11,3){\dashbox{1}(0,24)} 
	%\put(-7,0){$*$}
	\put(75,3){\circle*{6}}
	\put(59,19){\circle*{6}}
	\put(11,3){\circle*{6}}
	%\put(11,24){\circle*{6}}
	\put(59,3){\circle*{6}}
	%\put(-4,3){\line(1,0){12}} 
	\put(61,3){\line(1,0){12}} 
	\put(59,5){\line(0,1){12}} 
	\put(75,5){\line(0,1){12}} 
	\put(61,19){\line(1,0){12}} 
	\put(72,16){$*$}
\end{picture}
$1 \ ...\stackrel{\textstyle \frac{n'-2}{2} \strut}{n'-3\strut}\ \stackrel {\textstyle \underline{1} \strut}{\frac{n'-2}{2}\strut}$.
\\ 
\\
\item[(7)]
$A+E_6$.\\ 
$D_m|_{\Gam'}$:
\begin{picture}(23,18)
	\put(11,0){$*$}
\end{picture} \underline{1}.\\
$Z_{min}|_{\Gam'}$:
\begin{picture}(23,18)
	\put(11,0){$*$}
\end{picture} \underline{1}.\\
\item[(8)]
$A+D_5'''$: \\
$D_m|_{\Gam'}$:\begin{picture}(82,18)
	\put(11,3){\circle*{6}}
	\put(11,3){\line(1,0){12}}
	\put(27,3){\circle*{6}}
	\put(27,3){\line(1,0){12}}
	\put(43,3){\circle*{6}}
	\put(43,3){\line(1,0){12}}
	\put(59,3){\circle*{6}}
	\put(59,3){\line(1,0){12}}
	\put(70,0){$*$}
\end{picture}	
$1 \ 1 \  1 \  1 \ \underline{1}$.\\
$Z_{min}|_{\Gam'}$:\begin{picture}(23,18)
	\put(11,0){$*$}
\end{picture} \underline{1}.
\end{itemize}
Comparing the $D_m|_{\Gam'}$ with $Z_{min}|_{\Gam'}$, we can get this proposition.
\end{proof}

\subsection{Arithmetic genus of a singularity in the essentially irreducible case} 

Let $i$ be a non-negative integer and $p_f(V,o)=p>0$. We have
$p_a(iY)-1=i(p_a(Y)-1)+\frac{i(i-1)}{2}Y^2=m(i(p-1) +\frac{i(i-1)}{2}Z^2)$, where $m$ denotes the length of the Yau sequence for $Z$. In the degree one case, we can get a {\color{black}lower} bound {\color{black}for $p_a(V,o)$} by $\underset{i}{\max}\{ p_a(iY)\}=p_a(pY)=\frac{p(p-1)m}{2}+1$.

\begin{lemma}[Lemma 3.2, \cite{Ko1}] \label{lem3.10}
	$p_a(V,o) \ge \frac{p(p-1)m}{2}+1$ holds for a normal surface singularity $(V,o)$ of degree one, where $p=p_f(V,o)$ and $m$ denotes the length of the Yau sequence for $Z$.
\end{lemma}

In fact, we have the equality sign holds when $Z$ is essentially irreducible.

\begin{theorem} [i.e. Theorem \ref{mt2}] \label{thm3.11}
	Let $(V,o)$ be a normal surface singularity of degree one with $p_f(V,o)>0$, $Z$ {\color{black}be} the fundamental cycle on the minimal resolution. Assume that $Z$ is essentially irreducible, then $p_a(V,o)=\frac{p(p-1)m}{2}+1$, where $p=p_f(V,o)$ and $m$ denotes the length of the Yau sequence for $Z$.
\end{theorem}
\begin{proof}
By Lemma 3.1 in \cite{Ko1}, we have {\color{black}that} $A_i=D_i-D_{i+1}$ is a $(-2)$-curve with $A_i \cdot D_i=-1$ for $1\le i <m$ and $D_m=Z_{min}$. Let $C$ be a cycle whose support is in $\pi^{-1}(o)$ such that $p_a(C)=p_a(V,o)$. Let $C=C_1+ \dots + C_n$ be a chain-connected component decomposition, where $C_i$ is a chain-connected cycle and $\mathcal{O}_{C_j}(-C_i)$ is nef for $i<j$.

Let $A \le Z$ be an irreducible component such that $A$ is not {\color{black}a} $(-2)$-curve. Since $Z$ is essentially irreducible, we have $p_a(C-C_i)=p_a(C)-p_a(C_i)-C_i \cdot (C_1+ \dots + C_{i-1}+C_{i+1}+\dots +C_n)+1 \ge p_a(C)-p_a(C_i) +1 \ge p_a(C)+1$ when $A \nleq C_i$. So, we can get that $A \le C_i$ for any $i$.

We prove $p_a(C) \le \frac{p(p-1)m}{2}+1$ by induction on the length $m$ of the Yau sequence. When {\color{black}$m=1$}, we have $Z=Z_{min}$ and $Z \cdot A=-1$. Since the coefficient of $A$ in $Z$ is $1$ and $A \le C_i \le Z$, we have $A \cdot C_i \le A \cdot Z =-1$. Furthermore, since $\mathcal{O}_{C_j}(-C_i)$ is nef, we have $C_i \cdot C_j=A \cdot C_j +(C_i-A)\cdot C_j \le A \cdot C_j \le -1$ for $i<j$. Then 
$$p_a(C)-1=\sum_{i=1}^{n}(p_a(C_i)-1)+\sum_{i<j}C_i \cdot C_j \le n(p-1)-\frac{n(n-1)}{2} \le \frac{p(p-1)}{2}.$$
When $m>1$, assume the inequality holds for $m-1$, let us proceed to show {\color{black}that} it is true for $m$. We consider all the chain-connected component{\color{black}s} $C_i$ such that $A_1 \le C_i$ and assume that $n_0$ is the number of these cycles. Since the coefficient of $A_1$ in $Z$ is $1$ and $C_i \le Z$, we have $A_1 \cdot C_i \le A_1 \cdot Z=-1$ {\color{black}when $A_1 \le C_i$}. Furthermore, since $\mathcal{O}_{C_j}(-C_i)$ is nef, we have $C_i \cdot C_j = A_1 \cdot C_j+(C_i-A_1)\cdot C_j \le -1$ {\color{black}when $A_1 \le C_i$ and $A_1 \le C_j$($i<j$)}. Then 
	$$
	\begin{aligned}
	p_a(C)-1 &=(p_a(C-\sum_{A_1 \le C_i}C_i)-1)+(p_a(\sum_{A_1 \le C_i}C_i)-1)+(C-\sum_{A_1 \le C_i}C_i) \cdot (\sum_{A_1 \le C_i}C_i) \\
	&\le (p_a(C-\sum_{A_1 \le C_i}C_i)-1)+\sum_{A_1 \le C_i}(p_a(C_i)-1)+\sum_{A_1 \le C_i,A_1 \le C_j,i<j}C_i \cdot C_j \\
	&\le (p_a(C-\sum_{A_1 \le C_i}C_i)-1)+n_0(p-1)-\frac{n_0(n_0-1)}{2} \\
	&\le (p_a(C-\sum_{A_1 \le C_i}C_i)-1)+\frac{p(p-1)}{2}
	.\end{aligned}
	$$
Notice that $A_1 \nleq C-\sum_{A_1 \le C_i}C_i$ and $D_2=Z-A_1$, we know that $\textup{supp}\ (C-\sum_{A_1 \le C_i}C_i) \subseteq \textup{supp}\ (D_2)$. 
Since $D_2$ is the fundamental cycle on its support and the length of the Yau sequence for $D_2$ is $m-1$, by the induction hypothesis, we have $p_a(C-\sum_{A_1 \le C_i}C_i) \le \frac{p(p-1)(m-1)}{2}+1$. It means that $p_a(C) \le p_a(C-\sum_{A_1 \le C_i}C_i)+\frac{p(p-1)}{2} \le \frac{p(p-1)m}{2}+1$.

It follows from Lemma \ref{lem3.10} that we have $p_a(V,o)=\frac{p(p-1)m}{2}+1$.
\end{proof}

We want to get a similar formula for the higher degree case: $$p_a(V,o)=p_a( ([\frac{p-1}{d}]+1)Y )=\frac{dm}{2}(\frac{2p-2}{d}-[\frac{p-1}{d}])([\frac{p-1}{d}]+1)+1,$$ 
where $d=-Z^2$ and $[ a ] :=\max \{ n \in \mathbb Z | n \le a \}$ for real number $a$(Gauss symbol). However this formula doesn't hold in {\color{black}the} general case. For example, if $m=1$, $Z \neq Z_{min}$ and $p \ge d+1$, then 
$$p_a(  ([\frac{p-1}{d}]+1)Y )=p_a(  ([\frac{p-1}{d}]+1)Z) < p_a(  [\frac{p-1}{d}]Z+Z_{min})\le p_a(V,o).$$
Notice that, {\color{black}the condition $Z^2=-1$ implies that $D_m=Z_{min}$,}   so we need the restrictive condition $D_m=Z_{min}$ when $d>1$.

\begin{lemma} \label{lem3.12}
	Let $(V,o)$ be a normal surface singularity of degree two or degree three, $Z$ is the fundamental cycle on the minimal resolution. Assume that $Z$ is essentially irreducible and $Z=Z_{min}$, then 
$$p_a(V,o)=p_a( ([\frac{p-1}{d}]+1)Z )=\frac{d}{2}(\frac{2p-2}{d}-[\frac{p-1}{d}])([\frac{p-1}{d}]+1)+1.$$
\end{lemma}
\begin{proof}
There exists a cycle $D$ such that $p_a(D)=p_a(V,o)$ and $\mathcal{O}_{Z}(-D)$ is nef since the  negative definiteness of the intersection matrix. Assume $A$ is the irreducible component $A \le Z$ such that $A$ is not {\color{black}a} $(-2)$-curve and $k$ is the coefficient of cycle $A$ of $Z$. Assume $ak+b$ is the coefficient of cycle $A$ of $D$, where $a,b \in \mathbb Z$ and $0\le b <k$. Then we claim that $aZ \le D' <(a+1)Z$. 

If $D \nleq (a+1)Z$, let $B=\max (D-(a+1)Z,0)>0$, where $\max (D_1,D_2):= \sum \max (n_i,m_i)E_i$ for $D_1=\sum n_i E_i$ and $D_2=\sum m_i E_i$, then we have $D-B<(a+1)Z$ by the definition of $B$. For any irreducible components $C \le B$, we have the coefficient of $C$ of $D-B$ is equal to the coefficient of $C$ of $(a+1)Z$. Therefore, $C \cdot (D-B) \le C \cdot (a+1)Z$. Notice that $Z$ is essentially irreducible and $A \nleq B$ and $Z=Z_{min}$, we have {\color{black}that} $C$ is {\color{black}a} $(-2)$-curve and $C\cdot Z=0$. It means that $B$ consist{\color{black}s} of $(-2)$-curve and $B \cdot (D-B) \le 0$. Then we have
$p_a(D-B)=p_a(D)-p_a(B)-B \cdot (D-B) +1 > p_a(D)$, contradicting the maximality of $p_a(D)$.

If $aZ \nleq D$, let $B=\max (aZ-D,0)>0$, then we have $aZ \le D+B$ by the definition of $B$. There exists an irreducible component $C \le B$ such that $C \cdot B <0$ by the negative definiteness of the intersection matrix. Since $C \le B$, we know that the coefficient of $C$ of $D+B$ is equal to the coefficient of $C$ of $aZ$. Therefore, $C \cdot (D+B) \ge C\cdot aZ =0$. But $C \cdot B<0$, contradicting that $\mathcal{O}_{Z}(-D)$ is nef. So we get $aZ \le D' <(a+1)Z$.

If $D\neq aZ$, assume $D=aZ+B$, where $0 < B <Z$. We claim that $p_a(B)-1 \le \frac{b}{k}(p_a(Z)-1)$. By {\color{black}the} Riemann-Roch theorem, $p_a(B)-1=\frac{1}{2}(B^2+B \cdot K)=\frac{1}{2}(B^2+bA \cdot K)$ and $p_a(Z)-1=\frac{1}{2}(Z^2+kA \cdot K)$. Hence the claim is equivalent to $B^2 \le \frac{b}{k}Z^2$. When degree two or degree three case, we know that $1 \le k \le 3$ by $Z^2=kA\cdot Z$. If $k=1$, the claim is obvious by $0\le b<k$. If $k=2$, it means that $d=2$ and $\frac{b}{k}Z^2 \ge -1$, the claim is obvious too. If $k=3$, it means that $d=3$ and $0 \le b \le 2$. The claim is obvious when $b \le 1$, so we only need to observe $b=2$ case. Since $B^2+2A \cdot K=2p_a(B)-2$, we know that $B^2$ is even number. Then $B^2 \le -2 =\frac{b}{k}Z^2$.

By $p_a(D)=p_a(aZ)+p_a(B)+aZ \cdot B -1 \ge p_a(aZ)$, we have $p_a(V,o)=p_a(D)=p_a(aZ)+p_a(B)-1 +ab Z \cdot A \le p_a(aZ)+ \frac{b}{k}(p_a(Z)-1+ak Z \cdot A) \le p_a(aZ)+ (p_a(Z)-1+ak Z\cdot A)=p_a((a+1)Z)$. So we know that $\underset{i}{\max}\{ p_a(iZ)\} \le p_a(V,o)=p_a(D) \le p_a (a+1)Z$, hence $$p_a(V,o)=\underset{i}{\max}\{ p_a(iZ)\}=p_a( ([\frac{p-1}{d}]+1)Z )=\frac{d}{2}(\frac{2p-2}{d}-[\frac{p-1}{d}])([\frac{p-1}{d}]+1)+1.$$
\end{proof}

When the singular point is of degree two, according to Proposition \ref{prop3.9}, we can generalized Lemma \ref{lem3.12} to $m>1$ case.

\begin{theorem}[i.e. Theorem \ref{mt3}]  \label{thm3.13}
Let $(V,o)$ be a normal surface singularity of degree two with $p_f(V,o)>0$, $Z$ is the fundamental cycle on the minimal resolution. Assume that $Z$ is essentially irreducible and $D_m=Z_{min}$, then 
$$p_a(V,o)=p_a( ([\frac{p-1}{2}]+1)Y )=m(p-1-[\frac{p-1}{2}])[\frac{p+1}{2}]+1=[\frac{p^2}{4}]m+1.$$
\end{theorem}
\begin{proof}
According to Proposition \ref{prop3.9}, we only {\color{black}need} to prove that the theorem holds in each case. There exists a cycle $C$ such that $p_a(C)=p_a(V,o)$ and $\mathcal{O}_{Z}(-D)$ is nef since the negative definiteness of the intersection matrix. Assume $A$ is the irreducible component $A \le Z$ such that $A$ is not {\color{black}a} $(-2)$-curve and $k$ is the coefficient of cycle $A$ of $Z$.

In case (1), let $C=C_1+ \dots + C_n$ be a chain-connected component decomposition, where $C_i$ is a chain-connected cycle and $\mathcal{O}_{C_j}(-C_i)$ is nef for $i<j$. As in the proof of Theorem \ref{thm3.11}, we have $A \le C_i$ for $1 \le i \le n$. Since $C_i$ is chain-connected, we have {\color{black}that} $\textup{supp}\ C_i$ is connected. With the notations as {\color{black}in} Theorem \ref{thm3.8}, we know the coefficient of cycle $B$ of $Z$ is $1$ for each cycle $B$ in $\Gam'$. We denote the number of the cycle $B$ that $B \le C_i$ and $B$ in the first $A_{n'}$(resp. the last $A_{n'}$) {\color{black}by} $a_i$(resp. $b_i$) for $1 \le i \le n$. Let $c_j$(resp. $d_j$) be the number of $i$ that $a_i=j$(resp. $b_i=j$) for $0 \le j \le n'$. Then we have $\overset{n'}{\underset{j=0}\sum}c_j=\overset{n'}{\underset{j=0}\sum}c_j=n$ and $p_a(C)-1=\sum_{i=1}^{n}(p_a(C_i)-1)+\sum_{i<j}C_i \cdot C_j \le n(p-1)-\overset{n'}{\underset{j=0}\sum}\frac{c_j^2-c_j+d_j^2-d_j}{2} \le np-{\underset{j=0}\sum}\frac{c_j^2+d_j^2}{2}.$ 
It is easy to see that $p_a(C)$ reaches its maximum value when $c_0=c_1=c_2=\dots=c_{n'}=d_0=d_1=\dots=d_{n'}=[\frac{p-1}{d}]$ and $m=n'+1$, so the theorem holds in case (1).

In case (2) or case (3), there exists {\color{black}a} unique cycle $B \le Z$ such that $B \cdot Z=-1$. We consider the coefficient of cycle $B$ of $C$, denote {\color{black}by} $2a+b$($a,b \in \mathbb Z, 0\le b \le 1$). As in the proof of Lemma \ref{lem3.12}, $aZ \le C$. If $b=0$, we have $\mathcal{O}_{C-aZ}(-Z)$ is numerically {\color{black}trivial}, hence, $p_a(C)-1=(p_a(C-aZ)-1)+a(p-1)-\frac{a^2-a}{2} \le (p_a(C-aZ)-1)+[\frac{p^2}{4}]$ and $\textup{supp}\ (C-aZ) \subset \textup{supp}\ (D_2)$. If $b=1$, there exists an end $D$ in $\Gam'$ that connected with $B$. Denote the coefficient of cycle $D$ of $C$ {\color{black}by} {\color{black}$c$}, we have $c=a$ and $D\cdot C=-1$ since $-1\le D\cdot C=-2c+(2a+b)\le 0$, hence, $\mathcal{O}_{C-aZ-B-D}(-Z)$ is numerically {\color{black}trivial} and $p_a(C)-1 \le p_a(C-B-D)-1=(p_a(C-aZ-B-D)-1)+a(p-1)-\frac{a^2-a}{2} \le (p_a(C-aZ-B-D)-1)+[\frac{p^2}{4}]$. By induction, the theorem holds in case (2), case (3).

In case (4), case (5), case (6), we have $m=2$. There exists {\color{black}a} unique cycle $B \le Z$ such that $B \cdot Z=-1$, denote the coefficient of cycle $B$ of $C$ {\color{black}by} $2a+b$($a,b \in \mathbb Z, 0\le b \le 1$). As in the proof of Lemma \ref{lem3.12}, $aZ \le C$. The theorem is trivial when $b=0$, so we only need to prove the theorem when $b=1$. We have $p_a(C)-1=(p_a(aZ)-1)+(p_a(C-aZ-B)-1)+(B\cdot(C-B)-1)$. Notice that $B \cdot (C-B)-1 \le 1$, $p_a(aZ)-1\le [\frac{p^2}{4}]$, and $p_a(C-aZ-B)-1\le [\frac{p^2}{4}]$, we have $p_a(C)-1 \le 2[\frac{p^2}{4}]+1$. Whatmore, $p_a(aZ)-1 = [\frac{p^2}{4}]$ means $a \ge [\frac{p}{2}]$, and $p_a(C-aZ-B)-1 = [\frac{p^2}{4}]$ means $[\frac{p}{2}]D_2 \le C-aZ-B \le [\frac{p+3}{2}]D_2$, then $B \cdot (C-B)-1=B\cdot aZ+B\cdot (C-B-aZ)-1<-[\frac{p}{2}]+[\frac{p+3}{2}]-1 \le 1$. So the theorem holds in case (4), case (5), case (6).

In case (7), we have $m=3$. There exists {\color{black}a} unique cycle $B \le Z$ such that $B \cdot Z=-1$, denote the coefficient of cycle $B$ of $C$ {\color{black}by} $2a+b$($a,b \in \mathbb Z, 0\le b \le 1$). Similarly, we have $p_a(C)=p_a(aZ)+p_a(C-aZ-bB)-1\le p_a(C-aZ-bB)+ [\frac{p^2}{4}]$ when $b=0$, and $p_a(C)=p_a(aZ)+p_a(C-aZ-bB)-1+(B\cdot(C-B)-1)$ when $b=1$. Since $D_2$ is the same as case (4) with $k'=0$ and $n'=5$, we have $p_a(C-aZ-bB) \le 2[\frac{p^2}{4}]+1$, and $p_a(C-aZ-bB) = 2[\frac{p^2}{4}]+1$ means $C-aZ-bB \le [\frac{p+3}{2}](D_2+D_3)$. Hence, similarly as case (4)-(6), we have the theorem holds in case (7).
\end{proof}

\end{document}